\documentclass[preprint]{imsart}
\usepackage{graphicx}
\usepackage{natbib}
\usepackage{delarray}
\usepackage{hyperref}
\usepackage{amsthm}
\usepackage{rotating}
\usepackage{url}
\usepackage{colortbl,color}
\usepackage{anysize}
\usepackage{multirow,amssymb,amsmath, amsfonts}
\usepackage{lscape, subcaption}
\usepackage{verbatim}

\usepackage{enumitem}
\usepackage{bm}
\usepackage[toc,page]{appendix}
\usepackage{mathabx}

\newtheorem{thm}{Theorem}

\newtheorem{lemma}{Lemma}

\def\boxit#1{\vbox{\hrule\hbox{\vrule\kern6pt\vbox{\kern6pt#1\kern6pt}\kern6pt\vrule}\hrule}}

\def\argmin{\mathop{\rm argmin}}

\allowdisplaybreaks

\newcommand{\real}{I\kern-0.37emR}
\definecolor{annotxt}{rgb}{0,0,1}
\def\argmin{\mathop{\rm argmin}}

\newcommand{\foology}{{\em foology}}
\newcommand{\foo}{{\em foo} \,}
\newcommand{\barr}{{\em bar} \,}

\newcommand{\minimize}{\mathop{\mathrm{minimize}}}

\makeatletter

\newcommand{\Rmnum}[1]{\expandafter\@slowromancap\romannumeral #1@}
\makeatother

\usepackage{natbib}
\begin{document}

\title{\bf Inferactive data analysis}
\runtitle{The philosophy of selective inference}
\author{Nan Bi, Jelena Markovic, Lucy Xia and Jonathan Taylor}

\begin{abstract}
We describe {\em inferactive data analysis}, so-named
to denote an interactive approach to data analysis with an emphasis on inference
after data analysis.
Our approach is 
a compromise between Tukey's exploratory (roughly speaking ``model free'') and confirmatory
data analysis (roughly speaking classical and ``model based''), also allowing
for Bayesian data analysis. We view this approach
as close in spirit to current practice of applied statisticians and data scientists
while allowing frequentist guarantees for results to be reported in the scientific literature, or Bayesian results where the data scientist may choose the statistical model (and hence the prior) after some initial exploratory analysis.

While this approach to data analysis does not cover every scenario, and every possible
algorithm data scientists may use, we see this as a useful step in concrete providing tools (with frequentist statistical
guarantees) for current data scientists.

The basis of inference we use is {\em selective inference} \citep{exact_lasso,optimal_inference}, in particular
its randomized form \cite{randomized_response}. The randomized framework, besides providing additional power and shorter confidence intervals, also provides explicit 
forms for relevant reference distributions (up to normalization) through
the {\em selective sampler} of \cite{selective_sampler}. The reference
distributions are constructed from a particular conditional distribution
formed from what we call a DAG-DAG -- a Data Analysis Generative DAG.
As sampling conditional distributions in DAGs is generally complex, the
selective sampler is crucial to any practical implementation of inferactive
data analysis.

Our principal goal
is in reviewing the recent developments in selective inference as well as
describing the general philosophy of selective inference.

\end{abstract}

\maketitle

\begin{quote}
The idea of a scientist, struck, as if by lightning with a question, is far from the truth. -- \cite{tukey_both}.
\end{quote}

\section{Introduction} \label{sec: introduction}

We begin by describing what we might call a typical 
day / week / quarter in an applied statistician or data scientist's working life.
In this work, we will typically use the term {\em data scientist} rather
than applied statistician, not only because such terms are {\em en vogue} at the moment,
but, as we will see, our approach to data analysis actually considers the applied statistician
a scientist, running ``experiments'' of their own. 

We think of our data scientists as somewhat subordinate to a scientist running
experiments that sample some process.
With a dataset in hand, data scientists are tasked with discovering structure in it, reporting
this structure to their scientist colleagues and ultimately the scientific literature.
To achieve this goal, a data scientist will  commonly do some exploratory
data analysis given by a sequence of queries (in silico function evaluations) based on their 
interests and knowledge existing statistical toolboxes. 
The goal of such queries is exploratory in nature, perhaps trying to discover interesting
structure or features that might be used in modelling further down the pipeline.

Sometimes, based on the results from previous queries, they may decide to consult from sources such 
as existing literature, and make additional queries to extract more information.  
Based on this information, hypotheses of interest are formed, often with corresponding
parameters or targets of interest. For these parameters, a data scientist is tasked with producing
a report for the scientist to publish in the scientific literature.

It is at this point in the data scientist's every day work that data analysis conflicts with classical statistical inference.
While it is natural that a data scientist will want to explore their data to find hypotheses of interest, classical
mathematical statistics relies on having specified a statistical model before observing the data.
This conflict is well-documented in the literature. For instance, \cite{diaconis1981} compares this application of classical
mathematical statistics in this context to a primitive ritual, while referring to exploratory data analysis as a form {\em magical thinking}. Diaconis also references \cite{leamer} who also has some
interesting descriptions of this conflict. As pointed out by \cite{leamer} this conflict essentially predates
statistics, and was recognized by Sherlock Holmes:

\begin{quote}
It is a capital mistake to theorize before one has data. Insensibly one begins to twist facts to suit theories, instead of theories to suit facts. -- \cite{doyle_adventures_1892}
\end{quote}

\cite{diaconis1981} ultimately leaves some room for magical thinking in mathematical statistics in the form of a ``working hypothesis.'' In this work, the concept of ``working hypothesis'' is directly related to the statistical model of Section \ref{sec:model}
our data scientist chooses based on what he has observed about the data.

This conflict in the way a data scientist approaches data analysis and
the way a mathematical statistician produces tools for inference has very important implications
for how data is used to further science.
It needs no acknowledgment that reporting $p$-values or confidence intervals based on classical mathematical statistics
have no guarantee if the same data is used to both generate the hypotheses of interest
as well as construct the report. \cite{breiman2001} called this the ``quiet scandal of statistics'', pointing
out that statisticians know that the guarantees they provide only apply without data snooping. On the flipside,
data scientists are not completely innocent: 
\cite{leamer} referred to the use of such methods in science as a fundamental (perhaps original) sin of a data scientist.

Recent work in selective inference has attempted to address
this scandal and provide data scientists with tools to construct reports with at least some form 
of statistical guarantees. Roughly speaking, these approaches can be described as either
{\em simultaneous inference}  \citep{posi}, in which a certain class of statistical functionals
is determined before looking at the data and coverage guarantees are constructed
simultaneously over the entire class, or {\em selective inference} \citep{optimal_inference,exact_lasso,tibshirani_exact_2016} which
conditions on the outcome of some model selection query to construct a reference distribution. 

In this work,
we take the conditional approach, i.e. selective inference, recognizing that many data scientists do not have a clear enumeration
of possible questions of interest before they see a typical data set. Indeed, one of the great
things about science is that data can (and should) have the ability to make a scientist (and hence a data scientist)
change their mind about the data generating process in question. In the context of inference, this change of mind may
manifest itself by changing the statistical model used to construct the relevant statistical results. This notion
of introducing a different statistical model based on observations about a dataset was introduced
in \cite{optimal_inference} as a {\em selected model} in which the selected model was described formally as if chosen by an algorithm. In practice, this model can and should be chosen in consultation with
the scientist who measured the data along with the results of the data scientist's query. As mentioned above, this statistical model is a formal
version of \cite{diaconis1981}'s working hypothesis. After $U$ (for user) has chosen a suitable model for the data, $U$'s working hypothesis  must
be adjusted to reflect the fact that he has used the data to discover
this hypothesis. This adjustment is formally done by conditioning on what $U$ has
observed before forming this model. This conditioning  results in a new statistical model, the {\em selective model} 
described in Section \ref{sec:selective:model}. 

Our ultimate goal is to  describe a new ``theory'' of data analysis that add to the existing theories of \cite{diaconis1981}: {\em inferactive data analysis}
which uses the conditional approach of selective inference as its basic building block. In a formal sense, this form of inference is just classical
mathematical inference applied to a selective model.
Formal justification for such results are not presented here, and can be found in
\cite{selective_sampler,randomized_response,selective_bootstrap,cmu_bootstrap, markovic2017adaptive}. Our goal here is to provide
the reader with a description of the general viewpoint and the main concepts needed
to carry out this program. 

\subsection{Example: HIV resistance data} \label{sec:introduction:example}

As an illustration of the ideas discussed here, we consider at a real dataset studied in \cite{rhee2006genotypic,selective_sampler}. The authors studied the genetic basis of drug resistance in HIV using markers of inhibitor mutations to predict a quantitative measurement of susceptibility to several antiretroviral drugs. The goal is to find the mutations that predict responses to drugs. In particular, we take the protease inhibitor subset of their data (``HIV dataset'') and we are interested in one specific drug, Lamivudine (3TC).  There are 633 cases and 91 different mutations occurring more than 10 times in the sample. 

A data scientist is tasked with building a regression model for this data. 
Let $U$  denote our data scientist. Our user $U$ decides to first investigate which mutations have the largest marginal effect
in determining resistance to 3TC. This query can be answered by marginal screening \citep{sure_2008,genovese_screening}. 
The mutations with a marginal $Z$-statistic with value greater than 2.5 are 
\begin{verbatim}
[P35I, P39A, P41L, P43Q, P67N, P74I, P74V, P83K, P118I, P184V, P200A, P208Y, P210W, 
P211K, P215Y, P219E, P228H].
\end{verbatim}
Having observed the most important marginal effects, $U$ decides to run LASSO \citep{tibs_lasso} with features given by these mutations
as well as their interactions to discover whether there are any important interactions between mutations. The user $U$ uses
LASSO with a fixed value of the regularization parameter based on theoretical considerations (c.f.~\cite{negahban_unified_2012}), where $U$ 
uses an estimate of $\sigma^2$ from the full model with 91 mutations. The resulting mutations are 
\begin{verbatim}
[P67N, P83K, P118I, P184V, P210W, P215Y, P228H, P67N:P83K, P67N:P184V, P67N:P211K, 
P83K:P184V, P184V:P210W, 'P184V:P215Y', P200A:P210W].
\end{verbatim}

Given that $U$ has observed these facts about the data, what report should they produce? Which features should $U$ use in the model reported? We argue that $U$, in conjunction with the scientist can choose which model to report -- this
is the selected model. Having fixed a model, how should $U$ produce confidence intervals and / or $p$-values?

\subsubsection{What is valid inference: a tale of two data scientists}

Of course, a different data scientist (whom we will call $U'$) may have made different decisions -- data analysis can be a highly subjective enterprise. Suppose
$U'$ had instead chosen to run LASSO on all 91 mutations in the first stage instead perhaps using a similar choice of 
regularization parameter. The resulting mutations are 
\begin{verbatim}
[P41L, P62V, P65R, P67N, P83K, P151M, P181C, P184V, P210W, P211K, P215F, P215Y].
\end{verbatim} 
In a second step, the data scientist again runs LASSO
using the most commonly co-occuring mutations. The results are 
\begin{verbatim}
[P62V, P65R, P67N, P83K, P151M, P181C, P184V, P210W, P211K, P215F, P215Y, P41L:P184V, 
P41L:P210W, P62V:P184V, P62V:P215Y, P67N:P184V, P67N:P211K, P83K:P184V, P181C:P184V, 
P181C:P211K, P184V:P210W, P184V:P215Y].
\end{verbatim}

We see that the ``important effects''
that $U$ and $U'$ have discovered are different. 
This is to be expected as they 
evaluated different functions on the data. Hence, if they were to 
create reports with $p$-values or confidence intervals using
the standard methods from linear regression
to publish in the scientific literature, they will
typically have different variables in them. This begs
the question: which of $U$ and $U'$ has produced a valid report?

The short answer most statisticians
would give is most likely neither, as both $U$ and $U'$ have 
ignored the effect of selection. Selective inference provides tools that
would allow $U$ and $U'$ to account for selection. Are both reports now correct? In general, no. Confidence intervals and $p$-values are always
defined relative to some statistical model -- if the statistical model
is badly misspecified then even adjusted for selection it is likely
badly misspecified. On the other hand, suppose
that both $U$ and $U'$ will form intervals or test hypotheses about 
population parameters under the assumption that our 633 cases
were samplied IID from some population of HIV+ patients. 
In this case, as long as they have properly
accounted for selection both reports will be statistically valid in the large sample limit,
though readers of the scientific literature should be careful about
interpreting these parameters as if some underlying parametric model is correct. This caveat of course stands whether or not $U$ and $U'$ have used
the data to choose what to report.

\subsubsection{Is all exploratory analysis the same?}

While our framework allows for production of valid results 
for both $U$ and $U'$ (at least in certain situations), it
does not mean that we should be completely agnostic about the role
of the data scientist in producing these results.
For example, suppose $U$ is research scientist with many years
devoted to understanding and modelling resistance to different
drugs and $U'$ is a freshman student in a data science class.
Imagine then a data set for a new HIV drug for which the goal
is to understand resistance based on 
mutation pattern. Expert $U$ uses well-established methods
within the field to discover important mutations, while 
$U'$ runs through a list of different techniques
presented in class and others discovered online. 

After their exploration, suppose $U$ and $U'$ arrive at the 
same list of important mutations and decide that they would like
to report $p$-values in a given regression model with these mutations.
Whose results should we trust more? How is this reflected in the resulting report?
In the framework described below, we will require that both $U$ and $U'$ declare the exploration they have done and, from our description above,
$U$ will have a simpler description then $U'$. 

We might
represent this exploration as the sequence of in silico function evaluations
each has used yielding two different dependency graphs (understood here in the sense of computer science) with the data  as the top node. We see that $U$'s dependency graph will likely be simpler than $U'$'s.
Selective inference takes the viewpoint that we should condition on 
our exploration of the data (in this case the nodes reflecting functions evaluated
on the data) so as not to bias our final results. In this
context, we would expect to be conditioning on ``less'' in $U$'s graph
compared to $U'$. As conditioning on more decreases
{\em leftover information} \citep{optimal_inference}, we can expect this to result in $U$ having more powerful tests than $U'$. 

Is such a decrease in power reasonable and / or desirable? If we take the
viewpoint that one must condition after
exploration (so that we have consumed some  information), then this seems to be a reasonable outcome. 
We can therefore
attribute a cost, in terms of information, to exploratory analysis: if we are
really interested in inference, we should not squander information
unnecessarily. Viewing loss of information as a cost is one that
is certainly recognizable today: large scale internet companies owe
at least part of their financial success due to their skills in acquiring data from their many users. Even more clear is the cost
of doing science. Without grant funding, scientists would have much less data.

\subsection{Outline}

The rest of this paper is organized as follows. In Section  \ref{sec:conditioning} we discuss the role of conditioning in selective inference.
Section \ref{sec:model} discusses the role of a statistical model
in selective inference. Section \ref{sec:dag} combines the previous
two sections defining a DAG-DAG, the basic object used for constructing
relevant reference distributions. Section \ref{sec:randomization}
considers randomized selection algorithms with Section \ref{sec:selective:sampler}
describing how particular randomized versions of selective inference
allow for explicit descriptions of the appropriate reference distribution. Bayesian inference after selection is described in Section \ref{sec:selective:bayesian}, where the data scientist may choose the model, including the prior, after some initial exploration.
In Section \ref{sec:3TC}, we carry out a randomized version of
data scientist $U$'s analysis above, writing out the explicit reference
distribution. Some computational and theoretical details are given in Appendix \ref{sec:simple:problem} through an example of inference for a prototypical simple problem: inference
after thresholding a sample mean.

\section{The role of conditioning}
\label{sec:conditioning}

A reader new to the literature on selective inference may ask themselves why we should
condition on the results of queries about the data.
The short answer is that we recognize that humans are easily biased
by data \citep{false_positive,tversky}, while statistical methods strive to provide unbiased conclusions. Conditioning on what we have observed about the data, combined with traditional
methods of mathematical statistics allows us to provide unbiased
conclusions even after we have observed some functions of the data.

As illustration, we begin with perhaps the simplest statistical inference problem: testing a point null hypothesis. To
isolate the idea from any particular model, we describe the approach in a
fictitious scientific field {\em foology} (rhymes with zoology) in which
scientists are searching for a conjectured to exist \barr particle in a system
with many \foo particles.
If one prefers to think in more concrete terms, we might 
consider ourselves an internet company where \foo particles represent typical users and \barr particles
are big-spenders. Alternatively, we might imagine ourselves a cyber-security company where
\foo particles are normal users and \barr particles are intruders.

Assume now that data  $T \in {\cal T}$ (in the form of a single particle) has been collected in some experiment. Based on theoretical
considerations in the field of \foology, it is assumed that the law of $S$ when the sample is 
\foo is known and denoted by $F_0$. We call the statistical 
model $H_0=\{F_0\}$, the {\em null hypothesis of foology}.

The task 
of determining whether or not this sample is \foo or \barr is given to a data scientist, $U$.
We identify the data scientist as a random variable to acknowledge that different data scientists will analyze
data differently, though we do not consider this in our simple example.\footnote{In fact, allow for users to be randomized algorithms
running without human intervention. 
However, most of the time we condition on $U$ so treating it
as random is somewhat moot. If we really understood the user, we might  use this information.}

As this field of science is relatively new, no gold standard has emerged as the best test to distinguish
between \foo and \textit{bar}, though we assume that the data scientist has some collection $(S_i)_{i \in I}$ of test
statistics that have shown promise, with some having better power in different regions of ${\cal T}$.
At this point, $U$ may choose any one of the $S_i$'s, which we denote by
$s^U_0$ and carry out a hypothesis test by drawing replicates $(\bar{T}_b)_{b=1}^B$ 
from $F_0$ and comparing the empirical distribution of $s^U_0(\bar{T}_b)$ to the observed value of $s^U_0(T)$.
Assuming, without loss of generality that the $S_i$'s reject for large values, this produces a $p$-value (sending $B \to \infty$)
\begin{equation} \label{eq:p0}
	p_0(T) = p_0(T;s^U_0) =  F_0\left(\left\{t: s^U_0(t) \geq  s^U_0(T)\right\}\right).
\end{equation}
Such a $p$-value is certainly recognizable to all statisticians.

Alternatively, to get some guidance regarding which test statistic to use, $U$
might decide to query the data. 
By query, we mean she might compute $Q_1^U(T)$, and, based on the observed value
decide which of the $S_i$'s to use.\footnote{We use the subscript 1 in anticipation
of allowing a second query below. Below, we will also allow the query to involve possible additional randomization so
we might write $Q_1^U(T,\omega_1)$ where $\omega_1$ is drawn from some distribution known to $U$. }
In other words, our model of the  interactive aspect of data analysis proceeds in silico as a series of function evaluations where the
functions take $T$ as one of the arguments. The data scientist observe $q_1$, the return value of this function.

It is well documented that reporting a $p$-value such as \eqref{eq:p0}
with, $s^U_0$ replaced by the $S_i$ chosen after having observed $Q^U_1(T)$, is no longer appropriate as the test has been chosen based on the data. This is similar
to scenario \# 3 in \cite{gelman_forking}.
It is also reminiscent of the discussion in \cite{cox1958},
in which a statistician is presented with a draw from one of two possible normal distributions
with the same means but different variance. In this setting, the
statistician is also told which population the data is drawn from. Cox compares a conditional
test to an unconditional test and makes the point that if our objective is to determine ``what we can learn from the
data that we have'' then the conditional approach is clearly the right approach. 
A key distinction
between Cox's example and our data scientist is that the data scientist has decided to query
the data to acquire this ``additional'' data, while in Cox's example the
data scientist is simply told which population it was drawn from.

What is the analyst to report now? Having observed $Q_1^U(T)$ to be $q_1$, the distribution $F_0$ is
no longer an appropriate reference distribution:
the data analyst now knows that $T$ is in the event
$\left\{t:Q_1^U(t)=q_1\right\}$. 
A natural solution to this problem is the conditional approach \citep{cox1958}: the data analyst
can use whatever test statistic she chooses, say, $s^U_1$, as long as the reference distribution she uses 
is the restriction\footnote{A brief comment on notation and interpretation of the conditioning statements.
When conditioning on the return value of $Q^U_1(T)$ being $q_1$ we are really considering
the restriction of $F_0$ to this event. That is, we typically never construct the corresponding conditional laws
for any other values of $q_1$. Hence, throughout, we are implicitly using regular conditional probabilities
always evaluated at observed quantities.}
  of $F_0$ to $\{t:Q_1^U(t)=q_1\}$. This results in a $p$-value
\begin{equation} \label{eq:p1}
	p_1(T) = p_1(T;s^U_1,q_1) = F_0\left(\left\{t: s^U_1(t) \geq  s^U_1(T)\right\}\bigl|\:Q_1^U(T)=q_1\right).
\end{equation}
It is clear that under the null hypothesis of foology, i.e.~that $T$ is a \foo particle, then
such a $p$-value can be used to provide a Type I error guarantee. For instance,
assuming that $s^U_1$ has a density, then
\begin{equation} \label{eq:marginal:typeI}
	F_0\left(\left\{t:p_1(t) \leq \alpha\right\}\right) = \alpha.
\end{equation}
In fact, a stronger Type I error guarantee holds
\begin{equation*} \label{eq:selective:typeI}
	F_0\left(\left\{t:p_1(t) \leq \alpha \right\} \vert\: Q_1^U(T)=q_1\right) = \alpha.
\end{equation*}
We call this stronger Type I error guarantee a {\em selective Type I error} guarantee \citep{exact_lasso,optimal_inference}.

We can now summarize the {\em philosophy of selective inference} as:
\begin{quote}
Data scientists' choice of what to report
is easily biased by data snooping, invalidating
most guarantees of classical statistics.
Conditioning on the data scientist's observations allows data scientists
to produce unbiased statistical reports with similar (but selective) guarantees to reports produced by classical statistics if the classical methods
had been used properly.
\end{quote}

While this seems a natural solution to the problem, what has this solution provided us, and what has it cost us?
At its core, a hypothesis test is a model-based attempt to quantify the uncertainty
in $T$ to decide whether $T$ is a \barr particle or a \foo particle. Having observed
$Q^U_1(T)=q_1$ (because $U$ requested this function be evaluated) we have noted
that the variation in $T$ is modified -- generally speaking it has been diminished. 

Conditioning on the value of $Q^U_1$ returns the uncertainty to that determined by
$F_0$, in particular the uncertainty of $F_0$ restricted to the event $\{t:Q_1^U(t)=q_1\}$. 
As $F_0$ is known, this conditioning can be carried out exactly (in theory) and
the data scientist can avoid any implicit bias he or she may have introduced
based on observing $Q_1^U(T)$ is equal to $q_1$.
Hence, the ability to carry out ``unbiased'' frequentist inference without querying
the data (i.e.~knowledge of $F_0$) combined with conditioning on the result of the query allows us to carry out 
``unbiased'' frequentist inference after observing the result of the query. This seems a substantial gain
in that we recognize that data scientists rarely are ready to report a $p$-value without some exploration of the data.
While classical statisticians may want consumers of statistical methods to not explore the data at all
we recognize that this is generally not realistic. Our approach
 requires $U$ to declare the pair $(Q^U_1,q_1)$.
Of course,
in order to make this approach practical, we also must be able to construct the reference distribution 
in \eqref{eq:p1}. Constructing this reference in more realistic data analysis settings is considered in Section \ref{sec:selective:sampler} below.

We see then, that using $p$-values that control selective Type I error allows us 
to use the data to select which test statistic to use to test the null hypothesis of foology.
Let us contrast our $p$-value \eqref{eq:p1} with what we might call the {\em naive} $p$-value, which is referred to 
as ``Researcher degrees of freedom without fishing'' in \cite{gelman_forking}:
\begin{quote}
\dots computing a single test based on the data,
but in an environment where a different test would have been performed given different data;
thus $T(y; \phi(y))$, where the function $\phi(\cdot)$ is observed in the observed case.
\end{quote}
This $p$-value is
\begin{equation} \label{eq:p1:naive}
	p^{naive}_1(T) = p^{naive}_1(T;s^U_1,q) = F_0\left(\left\{t: s^U_1(t) \geq  s^U_1(T); Q^U_1 (t) = q \right\}\right).
\end{equation}
It should be clear that such a $p$-value does not provide any sort of {\em selective Type I error} guarantee.
In fact
\begin{equation*}
	p^{naive}_1(T) = F_0\left(\left\{t:Q_1^U(t)=q_1\right\}\right) \cdot p_1(T)
\end{equation*}
and hence
\begin{equation*}
	F_0\left(p^{naive}_1(T) \leq \alpha \:\vert\:Q_1^U(T)=q_1\right) = \frac{\alpha}{F_0\left(\{t:Q_1^U(t)=q_1\}\right)}.
\end{equation*}
The naive $p$-value can therefore be corrected here by dividing by an appropriate constant. We have therefore effectively
resolved this issue of ``researcher degrees of freedom without fishing".
Our corrected
$p$-value is just $p_1$ \eqref{eq:p1}.
We caution the reader that such a simple correction
will not generally work. In more complicated models naive $p$-values can similarly be defined though the relationship between
a selective $p$-value (i.e.~a $p$-value with selective Type I error guarantees) and a naive one is not
as simple as above. 

We saw above that the selective Type I error guarantee is in fact stronger than 
a marginal Type I error guarantee. As our $p$-value \eqref{eq:p1} provides control of this
error, by the no-free lunch principle, there must be some cost to using this $p$-value. This cost
comes in terms of power and can be seen when we revisit the idea that
$Q_1^U$ provides the data scientist ``additional'' information. 
In \eqref{eq:p1} we see that $Q_1^U$ is not really ``additional'' data -- it is a function of $T$. 
When  
conditioning on the observed value of $Q_1^U(T)$, the variation in $T$ is diminished.
This variation is of course the source of power of our statistical test, and is why \cite{cox1958}'s
 unconditional test has higher marginal power than the conditional test. This
loss of power can be framed in terms of the conditional information of $T|Q_1^U$ or the
{\em leftover Fisher information} \citep{optimal_inference} in parametric models.

\subsection{Do we have to worry about forking paths?}

Above, our user $U$ reported only $p_1(T)$, having chosen test statistic $s^U_1$ after observing
$Q_1^U(T)=q_1$. What would the user have done if they had observed $Q_1^U(T)=q_1'$? 
\cite{gelman_forking} referred to this problem as the {\em garden of forking paths}. 
What should we do to account for the possibility that $Q_1^U(T)$ {\em might have been $q_1'$ instead of $q_1$}?

One reply to this question is given in Section \ref{sec:explore:confirm} below, pointing out that the traditional exploratory / confirmatory model of data analysis never carry out such counterfactuals. 
Nevertheless, if we so desire, there is nothing
stopping us from going through such a counterfactual exercise, though
answering the question of what data scientist $U$ might have chosen if the
result had been $q_1'$ rather than $q_1$ is clearly a difficult problem.

Imagine then that $U$ has decided before hand what
test statistic she would use for every $q$ in the range of $Q^U_1$.
That is, suppose our data scientist $U$ describes a rule to choose a test statistic $S_{Q^U_1} : {\cal T} \times {\cal Q}^U_1$ as a function
where ${\cal Q}_1^U$ is the range of $Q^U_1$. 
This rule may be randomized or not. In theory, one might then construct a $p$-value, defined on all of ${\cal T}$ that provides a Type I error guarantee
like \eqref{eq:marginal:typeI} though the $p$-value itself would of course be a different
random variable. One way to form such a $p$-value is to simply report the
selective $p$-value for each possible outcome $q_1 \in {\cal Q}^U_1$, though in many scenarios, each $p$-value controls a different error rate making
their comparisons slightly difficult even though they can be computed in
theory. 

A second approach would be just to consider the 
marginal distribution of $S_{Q^U_1}(T, {\cal Q}^U_1)$ under $F_0$. In both cases we require our
data scientist to provide this map $S_{Q^U_1}$ which depends on the query they are going to evaluate. Description of this mapping is referred to as
{\em pre-registering} in \cite{gelman_forking}. It is also a required step in the simultaneous approach
of \cite{posi} in that they require knowledge of ${\cal Q}_1^U$ in order to describe the class of functionals
they seek to have simultaneous coverage guarantees.

Arguably, data analysis is quite subjective and
intuitive. Data scientists may not want to describe such a map $S_{Q^U_1}$. Indeed, it may be impossible for a user $U$ to describe their intuition in such formal terms. The conditional approach only
considers what the data scientist observes about the data, defining a $p$-value only on 
the {\em selection event}
\begin{equation*} \label{eq:selection:event}
	\left\{t:Q^U_1(t)=q_1\right\}.
\end{equation*}
In this sense, if we take the conditional approach we can ignore
the forking paths -- we are only interested in the value of queries for the ``data that we have'' as Cox might say.

If we decide to adopt this conditional approach, then, we are effectively pruning down the trees in this
garden of forking paths. While we are not requiring $U$ to describe all counterfactual data analyses
she might have done, we might be losing something. We have already acknowledged that we sacrifice some
power in this approach, though this is completely expected. We also are giving up on the possibility
of marginal guarantees. By this, we mean that without specifying a pre-registration map, there is no
sensible way to define a $p$-value on all of ${\cal T}$. 
However, as we argue in the following section, viewed in the large sense, the classical
confirmatory / exploratory paradigm has exactly the same issue.

Readers may point out that we are requiring $U$ to declare all of their
exploration as part of the process of reporting their $p$-value. Hence,
we are requiring $U$ to be honest. However, this is no different
than the simultaneous approach of \cite{posi}: data analysts can certainly try many transforms of features and covariates before forming
a design matrix $X$. The use of pre-registration
is a form of certification of honesty, which comes with at
least two costs. First of all, $U$
must specify a map such as $S_{Q^U_1}$. The second cost
is one we will discuss in Section \ref{sec:model} below when
we discuss the role of statistical models in our context. 
Using pre-registration, both $U$ and the scientist are prevented from
magical thinking in specifying a new statistical model {\em after}
having observed the results of $Q^U_1$. Exploratory
data analysis is meant to create tools that allows
humans to extract insight from data -- pre-registration
excludes the abuse of such information in specifying a statistical model, and hence
in producing reports about parameters that only became interesting
after some exploration of the data.

\subsection{The exploratory and confirmatory theory of data analysis}
\label{sec:explore:confirm}

Tukey argues that science invariably requires both exploratory data analysis: allowing $U$ to query the data and
observe the results; as well as confirmatory data analysis: reporting $p$-values and confidence intervals \citep{tukey_both,tukey_multiple}.
Of course, statisticians recognize the problems inherent in testing hypotheses using the same data
used to generate them. How then, can these conflicting goals be resolved?

The simplest approach
is to collect more data to evaluate the hypotheses generated in an exploratory phase. This is the classical exploratory / confirmatory theory of data analysis.

Let us suppose then that $U$ reports to the scientist collecting the data that test statistic $s^U_1$ seems like a promising
test statistic based on their exploratory data analysis, i.e.~the evaluation of $Q^U_1$ on $T$. The scientist
runs a second experiment in identical conditions to collect $T'$. As luck would have it, the science of
foology happens to tell us that $T'$  independent of $T$ and identically distributed. Hence, its distribution under
the null of foology is $F_0$. The data scientist then computes
\begin{equation} \label{eq:p0:conf}
	p'_1(T') = F_0\left(\left\{t': s^U_1(t') \geq s^U_1(T')\right\}\right).
\end{equation}
If the $p$-value is less than 0.05 -- they have discovered a \barr particle! ({\em at level 0.05}).

Such a $p$-value
is surely uncontroversial: as uncontroversial as any $p$-value or hypothesis test may ever be. We also
did not require $U$ to specify a map $S_{Q^U_1}$. 
However,
let us consider the totality of the uncertainty in this setting. In this setting, the
data we have sampled is $(T,T') \sim F \times F$ with the null hypothesis of foology expressed simply as $H_0:F=F_0$.
A little thought shows that
$p'_1$ is actually only well-defined on
\begin{equation*}
	\left\{(t,t'): Q^U_1(t)=q_1\right\}.
\end{equation*}
We say well-defined in the sense that the model used to construct a reference distribution here never
considers the possibility that $Q^U_1(T) \neq q_1$.
Even this uncontroversial exploratory analysis followed by confirmatory analysis does not consider
all the possible answers to the query evaluated on $T$. In this sense,
classical statistics never considers all of the
counterfactual forking paths of \cite{gelman_forking}.

A little more thought shows that our confirmatory $p$-value can be
expressed in terms of the law
\begin{equation*}
	{\cal L}_{F_0 \times F_0}\left(T'\:\big|\: T, Q^U_1(T)=q_1\right) \equiv {\cal L}_{F_0 \times F_0}\left(T'\:\vert\: T\right).
\end{equation*}
In the conditional approach we require only conditioning on the observed values of the queries, so it
is reasonable to construct a selective $p$-value testing the null hypothesis of foology using 
the law
\begin{equation} \label{eq:data:carve}
	{\cal L}_{F_0 \times F_0}\left((T',T) \:\vert\: Q^U_1(T)=q_1\right).
\end{equation}
In fact, tests constructed based on the law \eqref{eq:data:carve} can have significantly greater power than \eqref{eq:p0:conf}. This increase in power
is demonstrated numerically in \cite{optimal_inference} in a parallel data analysis paradigm: data splitting (c.f.~\citep{cox_split,hurvich_tsai}) in
which an exploratory and  confirmatory sample are produced from one larger sample.

\section{The role of the statistical model: Eureka!}
\label{sec:model}

Up to this point, our discussion in the field of foology has focused on detecting departures from the null hypothesis of
foology described by $F_0$. For scientists following the scientific method, the statistical model $H_0=\{F_0\}$ should be thought of as
a mathematical model of the current understanding of foology. Nothing in the scientific method says that this model is the
one model to rule them all. Indeed, one of the great things about science is that data can, and in some cases should, force
scientists to change their model. In this section, we discuss how inferactive data analysis addresses this issue.

For concreteness, suppose that based on $U$'s query, $B$ experiences a moment of magical thinking. Having observed $Q^U_1(T)$, $B$ posits a new
distribution $F_{E}$ ($E$ for Eureka!) as a possible competing model to $F_0$.
Using new data $T'$, $B$ asks $U$ to compute a new goodness-of-fit $p$-value 
\begin{equation} \label{eq:p0:eureka}
	F_E\left(\left\{t': s^U_1(t') \geq s^U_1(T')\right\}\right)
\end{equation}
though they are not restricted to the statistic $s^U_1$.
Of course, in this test, $B$'s goal is not to reject a null hypothesis as rejecting this null hypothesis
would tend to falsify the new model.

Some in the foology community are skeptical of $B$'s model. How are they to decide? Assuming the two distributions
$F_0, F_E$ have density with respect to a common measure, the natural construction uses the likelihood ratio by appealing to
Neyman-Pearson. Using this, $B$ and the rest  of the community can agree on a way to decide whether $T'$ has sufficient
evidence to falsify $F_0$ and conclude that $B$ has made a step forward in foology. This $p$-value would be
\begin{equation} \label{eq:p0:neyman}
	p^{classical} = F_0\left(\left\{t': \frac{f_E(t')}{f_0(t')} \leq \frac{f_E(T')}{f_0(T')}\right\}\right).
\end{equation}
Above, the terms $f_E, f_0$ are densities for $T'$ with respect to some common carrier measure. 

In constructing the above $p$-value, the scientists have essentially forgotten about the original data $T$.
In computing the $p$-value \eqref{eq:p0:eureka}, $B$ is asserting that
the law of $T'$ is $F_E$ under the new model. If $B$'s new theory can posit
a joint distribution for $(T,T')$ then a selective $p$-value using a construction similar to \eqref{eq:data:carve}.
As experimental consideration under the null hypothesis of foology indicated that $(T,T')$ are IID, it seems
possible that  the same holds under $B$'s new theory as well, though this is not necessary for the
construction of the selective $p$-value. Assuming existence of appropriate densities,
Neyman-Pearson indicates that the optimal test has the form
\begin{equation} \label{eq:p0:neyman:carve}
	p^{carve} = \frac{(F_0 \times F_0)\left(\left\{(t,t'): \frac{f_E(t,t')}{f_0(t,t')} \leq \frac{f_E(T,T')}{f_0(T,T')}, \: Q^U_1(t)=q_1\right\}\right)}
{(F_0 \times F_0)\left(\left\{(t,t'): Q^U_1(t)=q_1\right\}\right)}.
\end{equation}
(With some abuse of notation, 
we have used the same notation for the density of $(T,T')$ as the marginal density of $T$ or $T'$ above.)

In fact, $B$ need not even collect more data in order to construct a test with the same Type I error
guarantees as (\ref{eq:p0:neyman},~\ref{eq:p0:neyman:carve}). Using just $T$, $U$ can report
the $p$-value \eqref{eq:p1} and enjoy the same guarantees as 
the two $p$-values above.
\begin{equation} \label{eq:p0:neyman:selective}
	p^{select} = \frac{F_0\left(\left\{t: \frac{f_E(t)}{f_0(t)} \leq \frac{f_E(T)}{f_0(T)}, \: Q^U_1(t)=q_1\right\}\right)}
{F_0\left(\left\{t: Q^U_1(t)=q_1\right\}\right)}.
\end{equation}

The three $p$-values can be used to test exactly the same hypothesis comparing $F_0$ to $B$'s new theory $F_E$.
Thresholding each of the three $p$-value at $\alpha$ results in a test that controls selective Type I error.
We note that, \eqref{eq:p0:neyman:carve}
and \eqref{eq:p0:neyman} require the additional cost and time of collecting data $T'$.

In terms of power, Theorem 9 of \cite{optimal_inference} demonstrates
that in the context of data splitting \eqref{eq:p0:neyman:carve}
dominates \eqref{eq:p0:neyman} in quite general settings. In our case,
while $T$ and $T'$ were collected in 
different experiments, the sampling distribution of $(T,T')$ is equivalent
to one in which both are collected
contemporaneously and we use one sample to choose which test
statistic to use. 
We expect that \eqref{eq:p0:neyman:carve} dominates the other two while
\eqref{eq:p0:neyman} dominates \eqref{eq:p0:neyman:selective}. 

The test based on $p$-value \eqref{eq:p0:neyman:carve} or \eqref{eq:p0:neyman:selective}
violates the classical separation between confirmatory and exploratory data analysis as $B$ only
posited the model $F_E$ {\em after } $U$ had reported the results of query $Q^U_1$.
As in the simple point null case, conditioning
on what $U$ has observed permits $U$ to construct
statistical reports that are unbiased in the sense that the tests all enjoy control of the
selective Type I error.

\subsection{Selective model} \label{sec:selective:model}

Our example above considers a model $\{F_0, F_E\}$. In
more realistic examples, our statistical models
are typically more complex. Consider then a parametric model
for the law of $(T,T')$: 
$$
{\cal M}=\left\{F_{\theta}: \theta \in \Theta\right\}
$$
containing $H_0$.
Clearly, up to computational considerations, it is simple to
construct the Neyman-Pearson test of $H_0$ versus $H_A:F=F_{\theta}$
for any $\theta \in \Theta$. All of the usual theory of parametric
statistical models can be applied to the {\em selective model}
\begin{equation} \label{eq:selective:model}
	{\cal M}^*_{q_1} = \left\{F_{\theta|q_1}^*:\theta \in \Theta\right\},
\end{equation}
where 
\begin{equation*} \label{eq:selective:distribution}
	F_{\theta|q_1}^*(A) = F_{\theta}\left(A\:|\:Q^1_U(T)=q_1\right)
\end{equation*}
is the restriction of some $F_{\theta} \in {\cal M}$ to the selection event.

The model ${\cal M}^*_{q_1}$ is still a parametric model. Hence, any tools
available to us for parametric models are applicable in this selective
inference setting as well. As in our simple
example of testing the null hypothesis of foology, understanding
of ${\cal M}$ can be combined with conditioning on what we have observed about $T$ after $Q^U_1$ in order to construct {\em selectively valid} inference, i.e.~``unbiased'' inference after selection.

Of course, the construction of a selective model does not depend on 
${\cal M}$ being parametric as one can simply restrict each distribution in
some collection to the same event. 
Under appropriate conditions, notions such as consistency and weak convergence
along sequences of models transfer to corresponding sequences
of selective models. 
See \cite{optimal_inference,selective_bootstrap,randomized_response,cmu_bootstrap} for more complete descriptions of what structure of ${\cal M}$ transfers to ${\cal M}_{q_1}^*$.

\section{Putting the pieces together: the DAG-DAG} \label{sec:dag}

\begin{figure}[h!]
  \centering
a)    \includegraphics[width=0.4\textwidth]{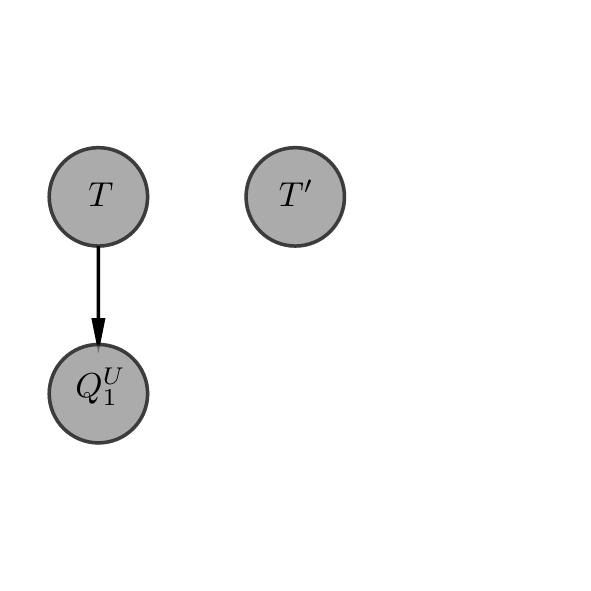}  \\
b)    \includegraphics[width=0.4\textwidth]{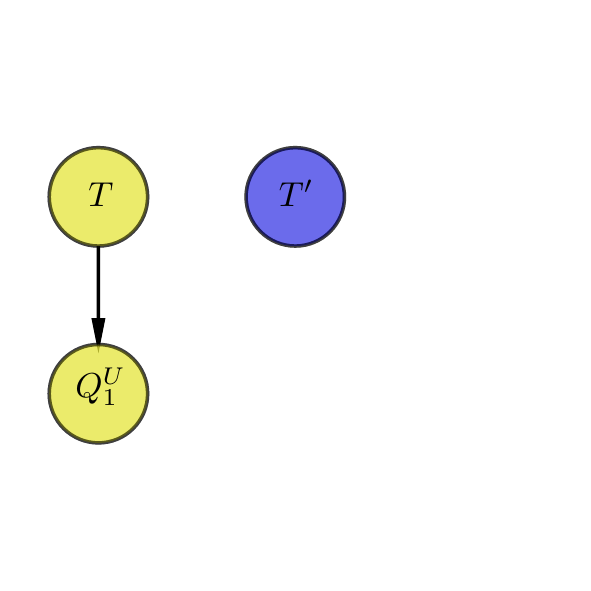} 
c)    \includegraphics[width=0.4\textwidth]{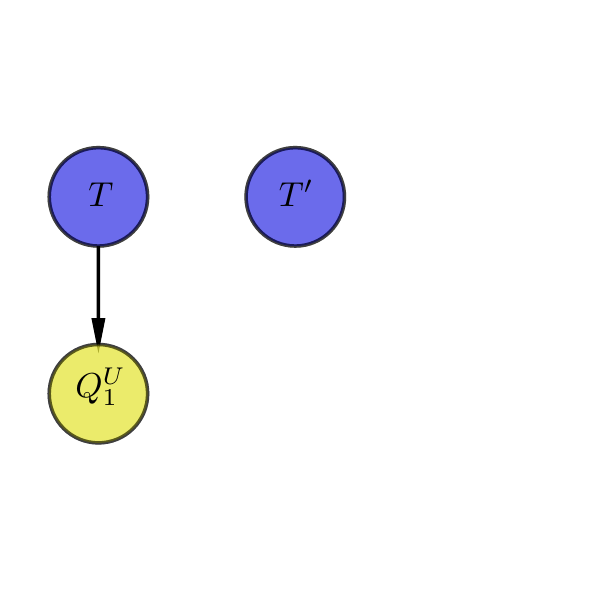}
    \caption{
a) The dependency graph for $(T,T',Q^U_1)$ in our foology example. Both the confirmatory \eqref{eq:p0:neyman} and carved $p$-values
\eqref{eq:p0:neyman:carve} have the same dependency graph.
b) The DAG-DAG for the confirmatory $p$-value \eqref{eq:p0:neyman}. The data analyst
posits a model for $(T,T')$, specifying
a model for $(T,T',Q^U_1)$ and conditions on the values of $(T,Q^U_1)$ (or, equivalently simply $T$).
c) The DAG-DAG for the selective $p$-value \eqref{eq:p0:neyman:carve}. The data analyst
posits a model for $(T,T')$, specifying
a model for $(T,T',Q^U_1)$ and conditions on the values of $Q^U_1$.
}
    \label{fig:dependency:graph}
\end{figure}

As the data is random, the results of functions evaluated (i.e.~the queries) on the data are also
random, hence any distribution we posit  for the data $T$ (such as $F_0$ or $F_E$) induces a joint distribution for $(T,T', Q^U_1)$. The random variables
 in this joint distribution are not independent -- $Q^U_1$ is a function of 
$T$. This dependence
can be expressed in terms of a {\em dependency graph} $G$
as in Figure \ref{fig:dependency:graph} a). 
Viewing the nodes of this DAG (Directed Acyclic Graph) as random
variables, we see that having specified a joint distribution for the data $(T,T')$, the
data scientist understands completely the joint distribution of $(T,T',Q^U_1)$. Specifying a statistical model ${\cal M}$ for $(T,T')$ (i.e.~a  collection of distributions) in turn specifies a statistical model for $(T,T',Q^U_1)$, the nodes
in the dependency graph $G$. This 
induced model is determined by ${\cal M}$ and the function $Q^U_1: {\cal T} \rightarrow {\cal Q}_1$. In terms of the dependency graph, we see then that 
the model is determined by a joint distribution for all observed data as
well as specification of the query functions which are represented by the
edges in the dependency graph. Let us call this induced model
${\cal M}' = {\cal M}'({\cal M}, Q^U_1)$.

As described in Section \ref{sec:conditioning}, in order to not bias herself
in her conclusions, the data scientist conditions on the value of $Q^U_1$.
This can be expressed as restricting each law in ${\cal M}'$ to the selection event, yielding a selective model ${\cal M}^* = {\cal M}^*({\cal M}, Q^U_1)$ as in \eqref{eq:selective:model}. Each distribution in ${\cal M}^*$ has the same dependency
structure (i.e.~DAG) as the dependency graph in Figure \ref{fig:dependency:graph} determined by $Q^U_1$ (which is encoded by the edge in the dependency graph) and the assumption of independence of $T$ and $T'$. 
This new model ${\cal M}^*$ represents the data scientists model of how
$(T,T',Q^U_1)$ are generated, conditional on what she has observed about the data. The null hypothesis of foology can be identified with one point in ${\cal M}^*$, and $B$'s Eureka model can be identified with another. 

We call such statistical models {\em DAG-DAGs}: Data Analysis Generative DAGs. They are specified by a dependency graph $G$ with vertices
of two types: data nodes and query nodes.
Along with $G$, the results of the queries are memorized as $Q$ and 
a statistical model ${\cal M}$ is posited for the data nodes. A
DAG-DAG then, is specified by the tuple $(G, Q, {\cal M})$.
The pair $(G,Q)$ is fully observed: it represents the queries
the data analyst makes about the data and their observed values. The statistical model ${\cal M}$ is specified by $U$ in order to carry
out statistical inference about distributions in ${\cal M}$. The
choice of ${\cal M}$ is allowed to depend on $(G,Q)$.
DAG-DAGs are the basic building blocks in describing inferactive data analysis and
are such that, up to sampling, they describe exactly how to reconstruct
the $p$-value \eqref{eq:p0:neyman:carve}. Figure \ref{fig:dependency:graph} b) describes this sampling distribution: the blue nodes indicate random variables for which the user must posit a joint distribution (coming from ${\cal M}$) 
while yellow nodes
indicate random variables whose observed values the data scientist conditions on (the values are stored in $Q$). Statistical model and observed query results
are integrated together through the dependency graph $G$.

 While control of the selective type I error requires the
data scientist to condition on which queries they have observed, the data
scientist may also condition on more (typically) at the cost of lower power \cite{optimal_inference}. Figure \ref{fig:dependency:graph} c) illustrates 
the sampling distribution used to construct \eqref{eq:p0:neyman}, constructed
by conditioning on the first experiment's data $T$.
Other reasons to condition on some of the data nodes 
in a DAG-DAG might be computational (c.f.~\cite{exact_lasso}) or related to other more classical statistical issues
in conditional inference.

In order to provide
valid selective inference, we will require the data scientist to make
note of the result of each query. The data scientist may request $B$
to collect additional data as the process evolves. In order to provide
valid selective inference, $U$ must specify a statistical model for the joint
distribution of any nodes in the DAG-DAG that are not observed.
One might
say that the data itself is ``observed'' and perhaps we should condition on it. However,
if we condition on the data itself, the resulting distributions in Figure \ref{fig:dependency:graph}
would all be degenerate, hence generating replicates from these distributions
would have no information about our null hypothesis.
Further, with the complexity of modern data sets, it seems fair to say that we never really
observe all of our data. Rather, we typically observe functions (often summary statistics) of the data.

In particular, if each query takes on
only countable many values, it is clear that the appropriate law at time $i$ is the
{\em selective distribution} $F^*_{0,i}, i \geq 1$ defined by the selective likelihood ratio \citep{randomized_response}
\begin{equation} \label{eq:likelihood:ratio:0}
	\frac{dF^*_{0,i}}{dF_{0,i}}(t^i) \:\propto\: \prod_{j=1}^i 1_{\{Q^{{\cal F}_j}(\cdot )=q^U_j\}}(t^i) \overset{\text{def}}{=} W_{(G_i,Q_i)}(t_i)
\end{equation}
where $q^U_j$ are the observed results of the query. In the above
likelihood ratio, the distribution $F_{0,i}$ is the statistical model
specified by the analyst while the right hand side can be read directly
off the dependency graph $G_i$ and observed query values $Q_i$ up to stage $i$. For a candidate
$F_{0,i}$ for the law of $T^i$, the actual selective likelihood ratio is
just the normalized $W_{G_i}$:
\begin{equation} \label{eq:likelihood:ratio:1}
	\frac{dF^*_{0,i}}{dF_{0,i}}(t^i) = \frac{W_{(G_i,Q_i)}(t^i)}{E_{F_{0,i}}[W_{(G_i,Q_i)}(t^i)]}
\end{equation}
Note also that the right hand side
depends only on observable quantities. Hence, \eqref{eq:likelihood:ratio:0}
defines a recipe to construct a selective model from the dependency graph and a model for $T^i$. For instance, after $B$'s Eureka moment, if
$U$ can posit a generative model $F^*_{E,i}$ for data collected up to query $i$, then she can construct
the selective distribution $F^*_{E,i}$ which has a density (with respect to $F_{E,i}$) proportional to \eqref{eq:likelihood:ratio:0}. More generally, 
given any statistical model for $T^i$, we can construct
a selective model from this dependency graph and the result of the queries
up to time $i$ using the right hand side \eqref{eq:likelihood:ratio:0}.

\subsection{Updating a DAG-DAG}

Having allowed a data scientist one query, of course
they will want to evaluate a second query, based on what they
have observed the value (i.e.~$Q^U_1$). Further, one data set does
not typically exist in a vacuum: labs will want to collect
more data based on the outcomes of earlier queries (including selective
hypothesis tests). If our framework is to be useful, it should 
behave nicely under such operations.

\begin{figure}[h!]
  \centering
a)    \includegraphics[width=0.4\textwidth]{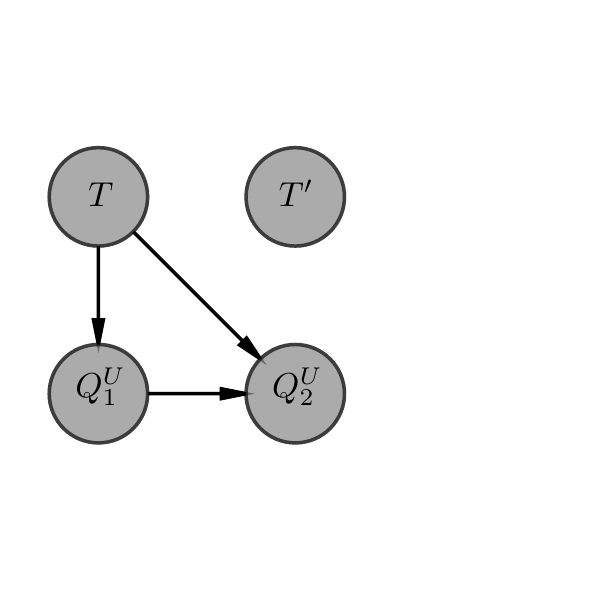}  
b)    \includegraphics[width=0.4\textwidth]{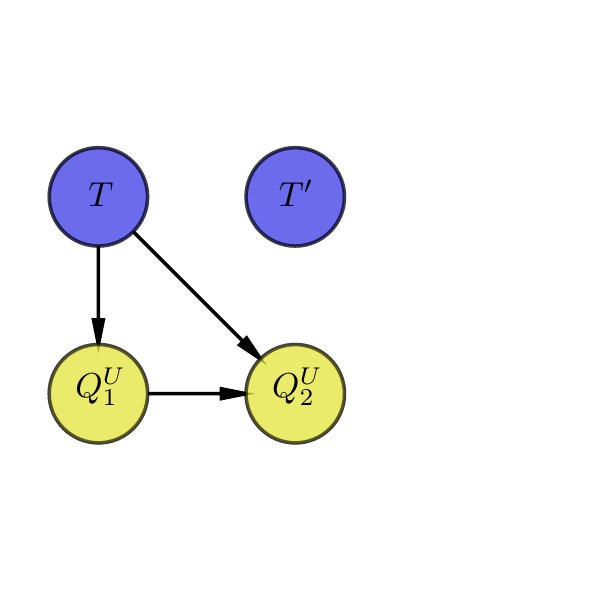}  \\
c)    \includegraphics[width=0.4\textwidth]{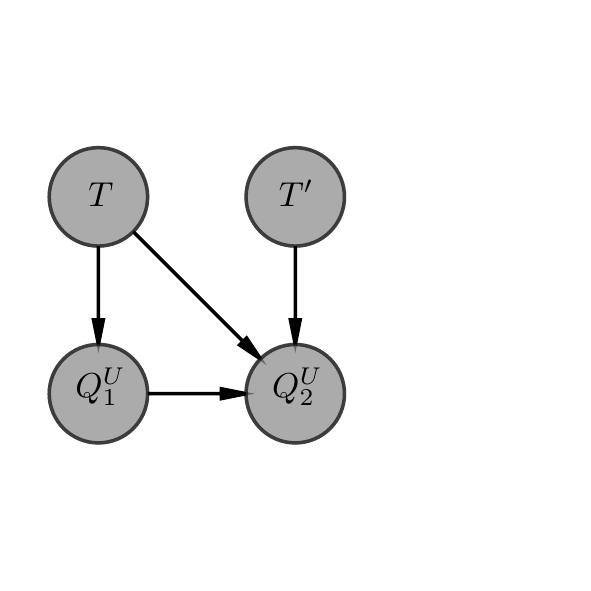}  
d)    \includegraphics[width=0.4\textwidth]{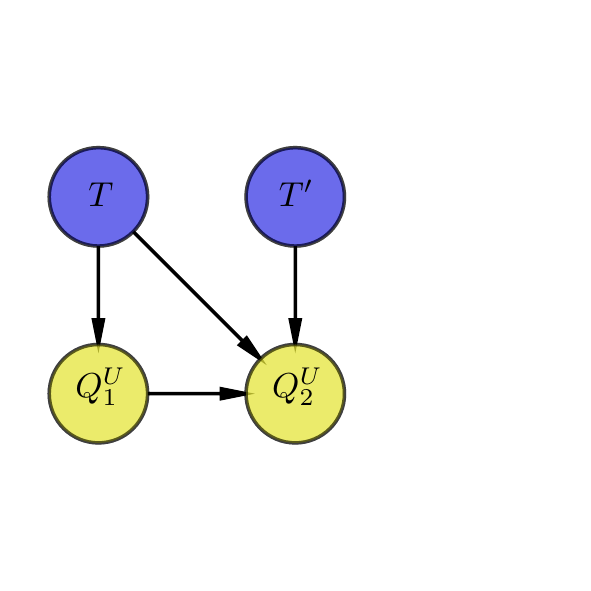}  
    \caption{
a) The dependency graph for $(T,Q^U_1,Q^U_2)$ after $U$ has evaluated
two queries on the data: $Q^U_2$ is a function only of $(T,Q^U_1)$.
b) The DAG-DAG used to compute the $p$-value $p_2$ \eqref{eq:p2}.
c) The dependency graph for $(T,Q^U_1,\bar{Q}^U_2)$ after $U$ has evaluated
two queries on the data: $\bar{Q}^U_2$ is a function of $(T,T',Q^U_1)$.
d) The DAG-DAG used to compute the $p$-value $\bar{p}_2$ \eqref{eq:p2:bar}.
}
    \label{fig:twoquery}
\end{figure}

In constructing \eqref{eq:p0:neyman:carve} we supposed that $B$
had commissioned collecting a second data set $T'$ in an effort
to decide whether $F_E$ is a true breakthrough in the field of foology.
In terms of the dependency graph, $Q^U_1$ is still a function
only of $T$, and we simply add an additional node
for $T'$. Specifying the DAG-DAG simply requires a joint
distribution for the pair $(T,T')$ which we have taken to be
IID $F_0$ under the null hypothesis of foology.
Figure 
\ref{fig:twoquery} a) and b) illustrate the dependency graph the DAG-DAG in which the data
scientist queries $T$ again, yielding $Q^U_2(T,Q^U_1)$. 

Of course, as $B$ has commissioned  $T'$  to be collected, $U$ might decide to evaluate a query
$\bar{Q}^U_2(T,T',Q^U_1)$. This would result in 
the dependency graph and DAG-DAG in Figure \ref{fig:twoquery} c) and d), respectively.
Having observed $Q^U_1$ and $Q^U_2$, the data scientist
chooses a test statistic and reports
\begin{equation}\label{eq:p2}
p_2(T,T')= p_2(T,T';s_2,[q_1,q_2]) = (F_0 \times F_0)\left(\left\{(t,t'): s^U_2(t,t') \geq  s^U_2(T,T')\:\bigl|\: Q^U_1(t)=q_1, Q^{U}_2(q_1,t)=q_2\right\} \right).
\end{equation}
Using query $\bar{Q}^U_2$, an appropriate  $p$-value might be
\begin{equation} \label{eq:p2:bar}
\bar{p}_2(T,T') = \bar{p}_2(T,T';s_2,[q_1,q_2]) = (F_0 \times F_0)\left(\left\{(t,t'): s^U_2(t,t') \geq  s^U_2(T,T') \:\bigl|\: Q^U_1(t)=q_1, \bar{Q}^{U}_2(t,t',q_1)=q_2\right\} \right).
\end{equation}

Clearly, this process could
continue indefinitely, collecting more data and recording
the result of $U$'s queries. 
Formally, then, we can 
think of this record of a data analysis
in a computer age version of the scientist's notebook in which a data scientist records their observations
at relevant steps of the data analysis.

Consider the expressions (\ref{eq:p0},~\ref{eq:p1},~\ref{eq:p2}): these are three 
different $p$-values for testing the null hypothesis
of foology. The first two
use only the data $T$ while \eqref{eq:p2} uses $(T,T')$. 
In what sense are these $p$-values, 
or hypothesis tests based on thresholding them
valid? Let $F$ denote the unknown true data generating mechanism for $T$. Clearly, 
thresholding $p_0$ at level $\alpha$ produces a valid test of $H_0:F=F_0$. It is also relatively simple
to see that similarly thresholding $p_1$ and $p_2$ at level $\alpha$ also result in 
valid tests of the above hypothesis. However, more is true. By construction, thresholding
$p_1$ or $p_2$ at level $\alpha$ produces a hypothesis test that is {\em selectively valid}
\citep{optimal_inference} at stages 1 and 2  of the data analysis, respectively.

By this we mean that at any stage $i$ in the data analysis, $U$ is able to produce a dependency graph $G_i$ describing all data collected up to this stage as well as all
function evaluations on this data. 
(Such practice is already encouraged in order to ensure in silico reproducibility of results, though such notebooks are not typically used for inference.)
Data scientist $U$, employed by $B$, may then report selective (for instance Neyman-Pearson) 
tests at any stage $i$ based on the graph $G_i$ and observed query values based on the law of their favorite test statistic (determined by $F_E$
if we use a Neyman-Pearson test) under the null hypothesis of foology.

We might denote the
$i$-th query by $Q^{{\cal F}_{i-1}}$ as it is measurable with respect to ${\cal F}_{i-1}$ and the data collected up to the $i$-th time point as $T^i$.
We denote the $i$-th test statistic by $s^{{\cal F}_{i}}$ as its selection is measurable with respect to ${\cal F}_{i}$. 
A natural choice for $i$-th $p$-value is just the usual one-sided $p$-value
with observed value $s^{{\cal F}_i}(T^i)$ and reference distribution based on the law
\begin{equation} \label{eq:conditional:family}
	s^{{\cal F}_i,E}(\bar{T}^i) \:|\: {\cal F}_i, \qquad \bar{T}^i \sim F_{0,i},
\end{equation}
where $F_{0,i}$ is the law of $\bar{T}^i$ under the null hypothesis of foology.

In this simple example, the statistical model used by $U$ is the two-point model
$\{F_0, F_E\}$ for each stage $i$ of the data analysis. In this sense,
we can make sense of both null and alternative hypotheses valid at each stage of the data analysis and $U$ can report a meaningful comparison of these
two hypotheses at any stage.

In more complex data analyses, it is 
likely that the model itself will change with $i$. For example, in our
3TC example, after marginal screening, $U$ might decide a reasonable model might include only variables that survive our thresholding rules, while $U$ might
change the statistical model if important interactions are discovered in the
second stage or if several variables are dropped in the LASSO analysis.
Nevertheless, $U$ can produce $p$-values after only marginal screening
or after the second LASSO analysis. In what sense are these valid?

The short answer is that the $p$-values at the $i$-th stage
are selectively valid given ${\cal F}_i$. In our foology
example, $B$ may refine the law $F_E$, allowing it to
depend on stage $i$ of the data analysis.
Alternatively, the field may have decided independently of $B$
that $F_0$ has been falsified
after $U$ has evaluated $Q^U_1$ but before producing $p_2$. The field of
foology may have then replaced $F_0$ with $\bar{F}_0$
a refined null hypothesis of foology. In this case, the appropriate
reference distribution for $p_2$ in the two-point model would seem to be
\begin{equation} \label{eq:selective:valid:0}
	{\cal L}\left(s_2(\bar{T}^2)  \:\bigl|\: {\cal F}_2\right), \qquad \bar{T}^2 = (T,T') \sim \bar{F}_{0} \times \bar{F}_0.
\end{equation}

In this case, $p$-values $(p_0, p_1)$ are selectively valid under $F_0$ while
$p_2$ is valid under $\bar{F}_0$. Such a situation is unavoidable if a field
of science is to refine its statistical model based on empirical data. The $p$-values
$(p_0,p_1)$ were indeed valid under the community's understanding of foology
at the time they were computed.

In a setting where
the dependency graph evolves as the field analyses data and collects more
data, a reasonable requirement of selective $p$-values are that they are {\em always-valid} \citep{johari2015always}.
That is,
\begin{equation} \label{eq:selective:valid}
	F_{0,i}\left(p_i \leq \alpha  \:\bigl|\: {\cal F}_i\right) \leq \alpha.
\end{equation}
However, even under $F_0$, the joint distribution of $(p_0,p_1,p_2)$ may be quite complex. For certain sequential model building data analysis 
schemes, the joint
distribution of such selective $p$-values can be described somewhat
concretely \citep{sequential_selective}.

\subsubsection{Merging DAG-DAGs}

It is certainly possible that more than one data scientist will analyze it.
Perhaps they want to work together on a joint analysis? Can we resolve
their two DAG-DAGs? Let $G$ be the dependency
graph of $U$ and $G'$ the dependency graph of $U'$. These graphs certainly share some nodes: the data nodes. They may also share other nodes: suppose
both $U$ and $U'$ use the same first query -- then (absent randomization
introduced below) both will have the node following data node $T^1$.
Allowing $U$ and $U'$ both access to $T^i$ at each stage $i$ is analogous to
broadcasting $T^i$ on a network where $U$ and $U'$ are different machines. Complete cooperation between $U$ and $U'$ then corresponds to broadcasting
all queries and observed values of the same network. 

More interesting questions certainly arise: what if $U$ is willing only to share part of their DAG-DAG with $U'$? What information must $U$ transmit to $U'$ to achieve this? What information must $U$ retain? Such questions are certainly interesting, though we leave these to further work.

\section{The role of randomization: experiments in data science}
\label{sec:randomization}

We have described how a data scientist can query a data set multiple
times before producing selectively valid $p$-values.\footnote{We have not yet explicitly described how to construct confidence intervals. These can of course be constructed by inverting selectively valid hypothesis tests. See
\cite{optimal_inference,randomized_response,selective_bootstrap,cmu_bootstrap} for further details.}
However, in the absence of collecting new data,
each query by the data scientist shrinks the support
of the relevant conditional distribution and it is not hard to imagine an inexperienced data scientist
querying the data so much that the conditional law begins to collapse towards a point mass. In this case,
the relevant hypothesis test will suffer a reduction in power to find \barr particles. In the selective
inference framework, this phenomenon has been  referred to as the ``leftover information'' as measured
by Fisher information in parametric selective models.

In the differential privacy literature \citep{reusable_holdout} and selective inference literature \citep{randomized_response}, a device that has been used to minimize the information loss
after each query is {\em randomization}. In this setting, instead of computing $Q^U_1(T)$ and $Q^U_2(T)$ the data scientist
computes a randomized version of these queries. That is, for the first query
the data scientist specifies a collection $(K^U_1(t))_{t \in T}$
of probability distributions on sample space $\Omega_1$ and makes a single draw $Q^U_1 \sim K^U_1(T)$. 
Based on
the value of $q^U_1$, the data analyst specifies a new family
of probability distributions $(K^{(U,[q_1])}_2(t))_{t \in T}$ and makes a single draw $Q^{(U,[q^U_1])}_2 \sim K^{(U,[q^U_1])}_2(T)$.

Simple examples of additional randomness include evaluating
functions on subsampled data, i.e.~data splitting methods \citep{hurvich_tsai,wasserman_roeder}. The area of adaptive data analysis (c.f.~\cite{reusable_holdout}) 
in computer science builds on randomization methods of differential privacy to provide some statistical guarantees in this interactive
setting.

The results of the queries here are random beyond the randomness
in the data $T$ or $(T,T')$. In this sense,
evaluation of the queries can be thought of as
experiments designed by $U$ to learn something about the
data. This can be a sticking point for many statisticians, as the results of the data analysis depend somewhat on the random seed
used to draw the queries. While we acknowledge this, we point out that
the results of the data analysis were already random to begin with: the data itself is random. Further, in a practical setting, 
$U$'s choice of queries (even if no additional randomness is inserted) is somewhat random and subject to $U$'s train of thought during the data analysis. These arguments will surely not convince every statistician of the usefulness of randomization. We present some more concrete advantages to randomization below: 
namely an increase in selective power, and, in many cases, a tractable
selective model.

In practical terms at the inference stage, the main difference is that the selective
likelihood ratio has been replaced by a (typically) smoother function.
Denote the kernel
at the $i$-th time point as $K^{{\cal F}_i}(t,\cdot)$ the selective likelihood ratio \eqref{eq:likelihood:ratio:0} is replaced by
\begin{equation} \label{eq:likelihood:ratio}
	\frac{dF^*_{0,i}}{dF_{0,i}}(t^i) \:\propto\: \prod_{j=1}^i K^{{\cal F}_j}(q^U_j;t^j).
\end{equation}
Computation of this likelihood ratio is generally impossible analytically, hence Monte Carlo methods are generally, though not always, necessary.
\cite{selective_sampler}, discussed below, presents some examples of randomization in which
the Monte Carlo burden can be substantially decreased, if not removed completely.

A common special case of the above setup is when the data scientist, based on observed information up to ${\cal F}_i$
draws $\omega_i$ from some distribution $G_i$ conditional on $T^i$ and the information up to ${\cal F}_i$. The data
scientist then observes $Q^{{\cal F}_j}=
h^{{\cal F}_j}(\omega_j)$, in which case
$$
K^{{\cal F}_j}(q^U_j;t^j) = G_j\left(\left\{\omega_j: h^{{\cal F}_j}(\omega_j)=q^U_j\right\} \:\bigl|\: t^j\right).
$$
At any stage $i$, the joint distribution $\bar{F}_{0,i}$ of
$(\bar{T}^i, \omega_1, \dots, \omega_i)$ (conditional on $U$) 
can be expressed in a simple form under $F_{0,i}$ via the likelihood ratio
\begin{equation}
\label{eq:likelihood:ratio:prob}
\begin{aligned}
\frac{dF^*_{0,i}}{dF_{0,i}}(t^i) &\:\propto\:\bar{F}_{0,i}\left(\left\{(t^i, \omega_1, \dots, \omega_i): h^{{\cal F}_j}(\omega_j)=q_j^U, 1 \leq j \leq i \right\} \right)\\
&\qquad = \prod_{j=1}^i G_j\left(\left\{\omega_j: h^{{\cal F}_j}(\omega_j)=q^U_j\right\} \bigl| t^j\right).
\end{aligned}
\end{equation}
The above form for the likelihood ratio means that draws from $F^*_{0,i}$ can be achieved by sampling from the appropriate marginal
of $F_{0,i}$ restricted to the event appearing on the right hand side.

For such randomized queries, in terms of the DAG-DAG only the functional form of arrows change. From deterministic functions of $T$, they are
now generated from kernels specified by $U$.

\subsection{Why randomize? Power and consistency}

The reader may be wondering why $U$ want to randomize our query? Will
we not decrease the chances of a {\em Eureka!} moment? 
While this certainly seems true, there is a tradeoff involved.
It has been empirically noted and demonstrated theoretically
in \cite{optimal_inference,randomized_response,selective_bootstrap}
that randomized selective inference 
procedures often enjoy benefits that unrandomized procedures do not have.

For example, under some sequences of selection procedures, it can be demonstrated that randomized algorithms remain consistent for much larger
classes of distributions than unrandomized procedures.
For instance, in the context of consistency and weak convergence,
the case of a single query was studied in detail in \cite{randomized_response}, including conditions for consistency
and a central limit theorem when
considering sequences of DAG-DAGs such that the randomization also depends on some sample size $n$. 
The regression examples presented in \cite{optimal_inference,randomized_response} also clearly demonstrate an increase in selective power after even
a small amount of randomization.

\subsection{Why randomize? The selective sampler}
\label{sec:selective:sampler}

We have seen above that our target of inference can be expressed
as a conditional distribution determined by a dependency graph $G$
and a statistical model for the relevant data in $G$. Sampling
from such conditional distributions is generally recognized to be a difficult
problem.
Perhaps we might be saved by special structure in our problems.
For instance,
\cite{exact_lasso, loftus_significance_2014-1, exact_screening} have shown that many common model selection procedures such as LASSO, forward stepwise and marginal screening have the selection event $\{Q^U_1 (t) = q_1\}$ expressed as explicit, if somewhat complicated, polytopes. \cite{optimal_inference} points out in the exponential family setting we further condition on the sufficient statistics of nuisance parameters to get an Uniformly Most Powerful Unbiased (UMPU) test. With randomization we need to evaluate the kernel $K^{{\cal F}_i}(t,\cdot)$ in the selective likelihood ratio in \eqref{eq:likelihood:ratio}. 
All these factors contribute to the difficulty of direct computation of the selective $p$-values. 

However, for many common queries, it turns out that the distributions
in the DAG-DAG can be reparametrized in a fashion that permits
explicit representation of relevant selective likelihood ratios
in terms of the data $T$ and certain auxiliary variables related to the query.
This is the subject of \cite{selective_sampler}, which we summarize here.

Consider a general randomized convex statistical program
\begin{equation} \label{eq:optimization:problem:rand}
	\hat{\beta}(T,\omega) = \argmin_{\beta \in \real^p}\quad \ell(\beta;T) + {\cal P}(\beta) -\langle \omega, \beta \rangle + \frac{\epsilon}{2} \|\beta\|^2_2,
\end{equation}
where $\ell$ is some smooth loss involving the data, ${\cal P}$ is some
structure inducing convex function,\footnote{The two most common examples of interest in statistical learning are 
${\cal P}(\beta) = h_K(\beta) = \sup_{\nu \in K} \nu^T\beta$ 
for some convex $K \ni 0$, i.e.~a seminorm, and that
${\cal P}(\beta) = I_{K}(\beta) = \begin{cases} 0 & \beta \in K \\
  \infty & \beta \not \in K \end{cases} $,
where $K$ is some cone or of the form $\left\{b:\|b\| \leq 1\right\}$ for some seminorm $\|\cdot\|$.
They address most of the statistical programs we are interested in.}
$\epsilon > 0$ is small that is sometimes necessary to assure the program has a solution and $\omega \sim G$ is a randomization chosen by the data analyst. Our query is represented as some function
related to the solution to \eqref{eq:optimization:problem:rand}, for instance
the sparsity pattern 
\begin{equation} \label{eq:selection:lasso}
	Q^U_1(T,\omega) = \left\{(j, \text{sign}(\hat{\beta}(T,\omega)_j): \hat{\beta}(T,\omega)_j \neq 0 \right\}
\end{equation}
if our penalty is a sparsity
inducing penalty such as the $\ell_1$ penalty in the LASSO \citep{lasso}. 
Our goal is to sample
\begin{equation*} \label{eq:inference:problem:rand}
	(T,\omega) \:|\: Q^U_1(T,\omega) = q_1.
\end{equation*}
Marginalizing over $\omega$ yields the desired selective distribution \eqref{eq:likelihood:ratio:prob}.

Of course, it is well known that the event \eqref{eq:selection:lasso}
can be expressed in terms of the KKT conditions of \eqref{eq:optimization:problem:rand}.
Inspection of the KKT conditions or subgradient equation
of \eqref{eq:optimization:problem:rand} yields the {\em reconstruction map}
\begin{equation*}
\omega = \phi(\beta,\alpha,z)  = \hat{\alpha}(s,\omega) + \hat{z}(t,\omega) + \epsilon \cdot \hat{\beta}(t,\omega)
\end{equation*}
defined on domain  
$$
\{(\beta, \alpha,z)\in {\cal B}_{q_1}(t):  \alpha + z + \epsilon \cdot \beta \in \text{supp}(G) \},
$$
where ${\cal B}_{q_1}(t)$ is defined by the selection event. In the
case of the LASSO with parameter $\lambda$ in which we write $q_1=(E,s_E)$ in terms of
variables and signs
$$
{\cal B}_{q_1}(t) = \left\{(\beta, \alpha, z): \beta_{-E}=0, 
\text{sign}(\beta_E)=s_E, \|z_{-E}\|_{\infty} \leq \lambda \right\}
$$
which does not even depend on $t$, though for other programs it may.

The subgradient equation gives rise to a canonical map
\begin{equation} \label{eq:canonical:map}
	\psi : (t,\omega) \rightarrow (t, \beta,\alpha,z)
\end{equation}
which yields a change-of-variable formula to sample the augmented parameter space $(t, \beta, \alpha, z)$ instead of the original $(t, \omega)$. The advantage is that ${\cal B}(t)$ usually has simple constraints on the optimization variables, e.g.~quadrant for $\beta$ and box for $z$, that are easy to sample from instead of the complicated constraints on the original problem. Further, choosing $G$ to have a Lebesgue density supported on all of $\real^p$ typically 
ensures the support is simple.

We can now state the main theorem from \cite{selective_sampler}.

\begin{thm}[Selective sampler]
\label{thm:change:measure}
The selective law
\begin{equation*} \label{eq:change:measure}
	{\cal L}_{F \times G}((T,\omega) \:|\: \hat{\beta}(T,\omega) \in A ) = {\cal L}\left(\left(T, \epsilon \cdot \beta +\alpha+z\right) \:|\: (\beta,\alpha,z)\in {\cal B}(T)\right)
\end{equation*}
for suitable ${\cal B}(T)$.
With $f$ denoting the density of $T$ and 
$g(\cdot|t)$ denoting the density of $\omega$ given $t$, the random vector
$(T,\beta,\alpha,z)$ has density proportional to
\begin{equation*} \label{eq:change:measure:density}
	f(t) \cdot g(\epsilon \cdot \beta +\alpha + z|t) \cdot \left|J\psi(t,\beta,\alpha,z)\right| \cdot 1_{\mathcal{B}(t)}(\beta,\alpha,z) \\
	=f(t) \cdot g(\phi(\beta,\alpha,z)|t)\cdot  \left|\det(D_{(\beta,\alpha,z)} \phi)\right| \cdot 1_{\mathcal{B}(t)}(\beta,\alpha,z) 
\end{equation*}
with the Jacobian the derivative of the map $\psi$ with respect to $(\beta,\alpha,z)$ on the fiber over $t$ 
\begin{equation*}
	\left\{(\beta, \alpha, z): Q^U_1(t,\epsilon \cdot \beta +\alpha + z)=q_1 \right\}.
\end{equation*}.
\end{thm}

\cite{selective_sampler} has abundant examples varying both the penalty as well as the loss function, including the canonical LASSO, forward stepwise and stagewise algorithms, marginal screening and generalized LASSO for the penalties, as well as squared-error, logistic and log-det for covariance matrix estimation for loss functions. 

\begin{figure}[h!]
  \centering
a)    \includegraphics[width=0.4\textwidth]{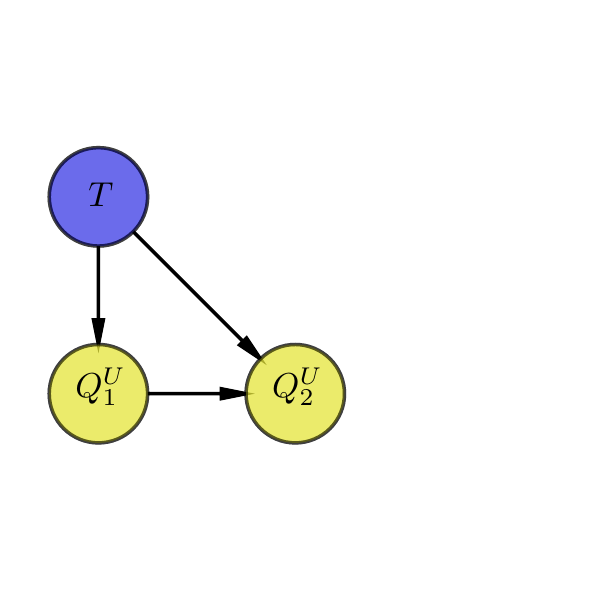}  
b)    \includegraphics[width=0.4\textwidth]{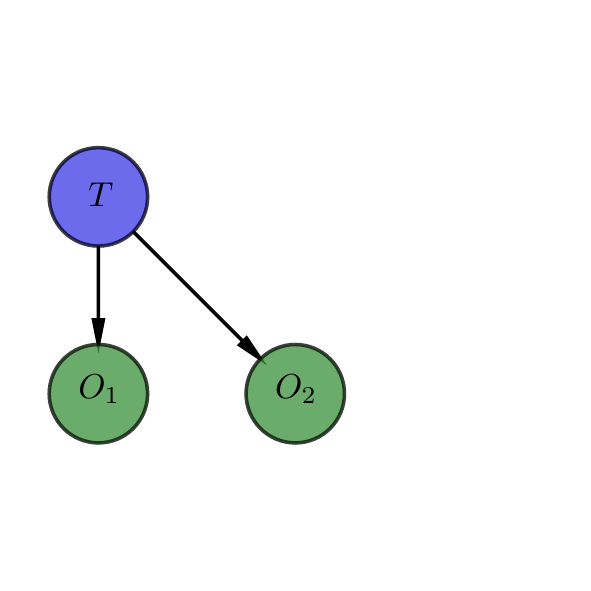} \\ 
c)    \includegraphics[width=0.4\textwidth]{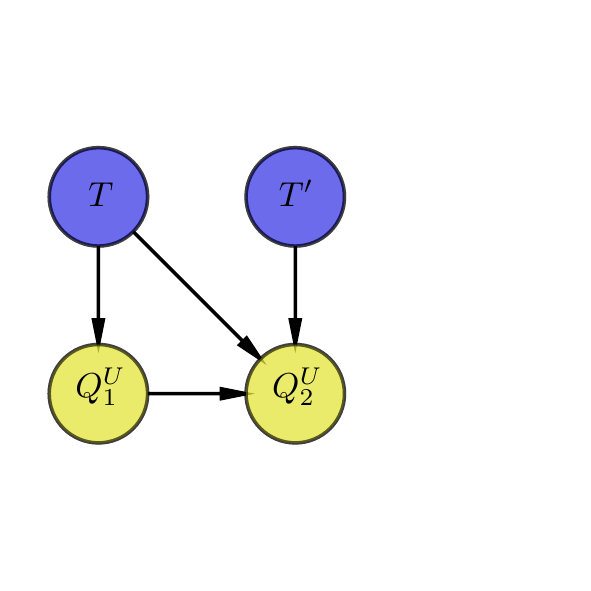}  
d)    \includegraphics[width=0.4\textwidth]{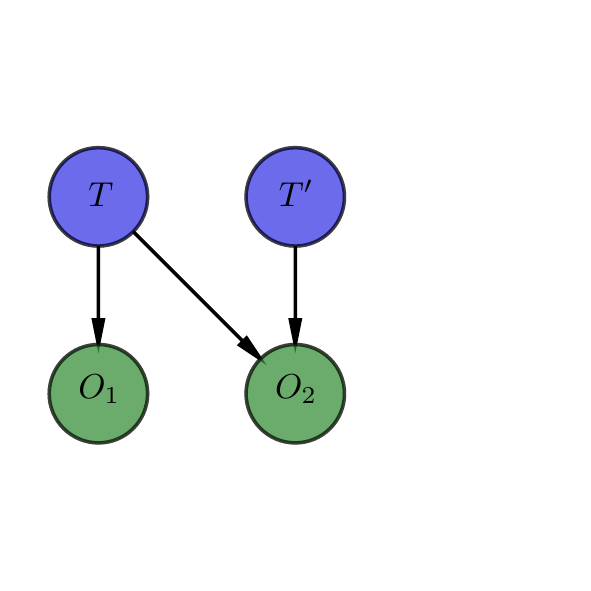}  
    \caption{
a) The DAG-DAG for computing a $p$-value
after evaluating two queries on $T$.
b) The selective sampler reparametrization of a) using
\eqref{eq:conditional:independence}.
c) The DAG-DAG for computing $p$-value \eqref{eq:p2:bar}. d) The
selective sampler reparametrization of c) using \eqref{eq:conditional:independence}.
}
    \label{fig:selective:sampler}
\end{figure}

We see then if each query is of this form, even if the objective
at stage $i$ depends on the queries at stages previous to $i$, the selective density is proportional to 
\begin{equation} \label{eq:conditional:independence}
	f_i(t^i) \cdot \left(\prod_{j=1}^i g_j\left(\epsilon_j \beta_j + \alpha_j + z_j \right) \cdot J\psi_j\left(t^j,\beta_j,\alpha_j,z_j\right) \right)
\end{equation}
supported on 
$$
\bigsqcup_{t^i \in \text{supp}(F_i)} \left(\prod_{j=1}^i \left\{(\beta_j, \alpha_j, z_j) \in {\cal B}_{q_j}(t^j): \beta_j \in \real^{p_j}, \alpha_j \in \partial \ell_j(\beta_j;t^j), z_j \in \partial {\cal P}_j(\beta_j) \right\} \right)
$$
with $f_i$ some point in a statistical model ${\cal M}_i$ provided by $U$ for the data collected up to stage $i$.
Note that the variables for each convex query are conditionally independent given $t^i$, hence can be sampled in parallel. This reparametrization 
is illustrated in 
Figure \ref{fig:selective:sampler} in which the appropriate reference distribution samples optimization variables $O_1=(\beta_1,\alpha_1,z_1)$ and 
$O_2=(\beta_2, \alpha_2, z_2)$ according to \eqref{eq:conditional:independence}.
The DAG-DAG on the left represents the sampling distribution when
$U$ evaluates the queries on $T$ without collecting new data. The DAG-DAG on the right contains green nodes which are marginalized over but do not require
$U$ to posit a distribution as their distribution is determined
by the law of $T$ and the selective sampler likelihood ratio
 \eqref{eq:conditional:independence}.

\subsection{A side benefit of randomization: counterfactual analyses}

Clearly, at any stage $i$, any other data scientist $U'$ can take $U$'s dependency graph $G_i$ and observed query results and specify their own statistical model. This is already done in the exploratory / confirmatory paradigm as well, and may often be requested in a peer-review process. In the context
of the DAG-DAG, this corresponds to simply changing the statistical
model posited for each of the data nodes in $G_i$. 

It is difficult, if not impossible, to have $U$ answer the question ``what query would you have chosen for $Q^U_2$ if the observed value of $Q^U_1$ were different?'' If $U$ were willing to pre-register such 
a randomized data analysis, the pre-registration declaration
effectively provides a mathematical answer to this question and $U'$ could draw the results of a second query having changed the result of query $Q^U_1$. Of course, $U'$ will have to declare how different or perhaps some specific value of $q_1$. In effect, this allows $U'$ to jump from one of Gelman's
forking paths to another, though the exercise seems to lead down to a rabbit hole if the parameter that $U$ reports a $p$-value about depends heavily on the answers to the previous queries. This rabbit hole is also avoided in the exploratory / confirmatory theory of data analysis because the field typically considers the exploratory data as fixed, i.e.~they have long since agreed to condition on it and not vary it.

Whether or not $U$ has preregistered their analysis, it is certainly possible to construct the reference distribution for different observed values of the queries (i.e.~fixing the functions that are evaluated and changing their return values). This makes it possible
for $U'$ to evaluate evidence for or against
other data generating mechanisms knowing only $G_i$. It is also certainly possible
for reviewers to take $U$'s statistical model and investigate how sensitive
inference is to the particular observed values of the queries.
Such counterfactual analyses are impossible without having randomized the query: for our observed $T$, either $Q^U_1(T)=q_1$ or it does not.



\section{The role of mathematical statistics} \label{sec:selective:bayesian}

\begin{quote}
	We comfortably divide ourselves into a celibate priesthood of statistical theorists, on the one hand, and a legion of inveterate data analysts-sinners on the other. -- \cite{leamer}
\end{quote}

In this work, we have tried to describe a way that allows data
scientists to explore their data to discover interesting hypotheses.
Conditioning on their path allows them to produce statistical reports
with selective guarantees. Formally,
the statistical theorist has told the data scientist
that the penance for their sins is that they should use the selective model
${\cal M}^*$ from their DAG-DAG for inference.

Having agreed that the proper reference distribution is found from the 
DAG-DAG, what role does the mathematical statistician have?
One prominent role is to further simplify the reference distributions
in the DAG-DAG. The selective sampler of Section \ref{sec:selective:sampler} simplifies
the conditional distributions of the DAG-DAG by explicit potentials
depending on the particular convex program solved (and which function of the solution the data scientist has observed). The mathematical statistician
can further simplify this when the queries are such that asymptotic
distribution theory can be applied.

\subsection{Selective Central Limit Theorem}
\label{sec:CLT}

To this point, we have considered fairly two-point statistical models in foology
without reference to any asymptotic regime. While we have
briefly defined parametric selective models, we have not gone
into much detail. More realistic
data analyses parametric models or nonparametric
models based on CLTs along sequences of models are likely to be important. 
In this section, we describe some of the asymptotic results found in \cite{randomized_response,selective_bootstrap}.

Suppose then that our data $T=T_n$ is such that, under statistical
model ${\cal M}_n$, the law of $T_n$ is asymptotically Gaussian. As our example below fits into this framework, for which we only observe
one data set (i.e.~$T$) we assume that we do not have access to $T'$ though in practice we of course may have such data. If the joint law of $(T_n,T'_n)$ were asymptotically Gaussian along the sequence of models ${\cal M}_n$ then much of what follows extends naturally. The data scientist, based on 
observing the outcomes of queries $Q^U_1, Q^U_2$ decides to report
a selective hypothesis about a parameter $\theta=\theta_n:{\cal M}_n \rightarrow \mathbb{R}$. In many contexts, $U$ will have access to an estimator
$\hat{\theta}_n$ of $\theta_n$ such that $n^{1/2}(\hat{\theta}_n - \theta_n(F))$ is
approximately centered Gaussian with variance 
$n \text{Var}_F(\hat{\theta}_n)$
for all $F \in {\cal M}_n$. 

\begin{figure}[h!]
  \centering
a)    \includegraphics[width=0.4\textwidth]{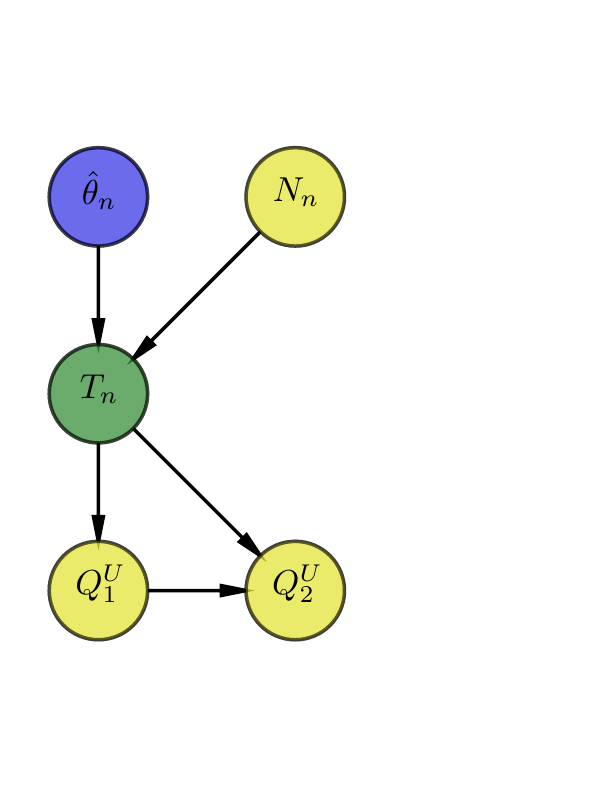}  
b)    \includegraphics[width=0.4\textwidth]{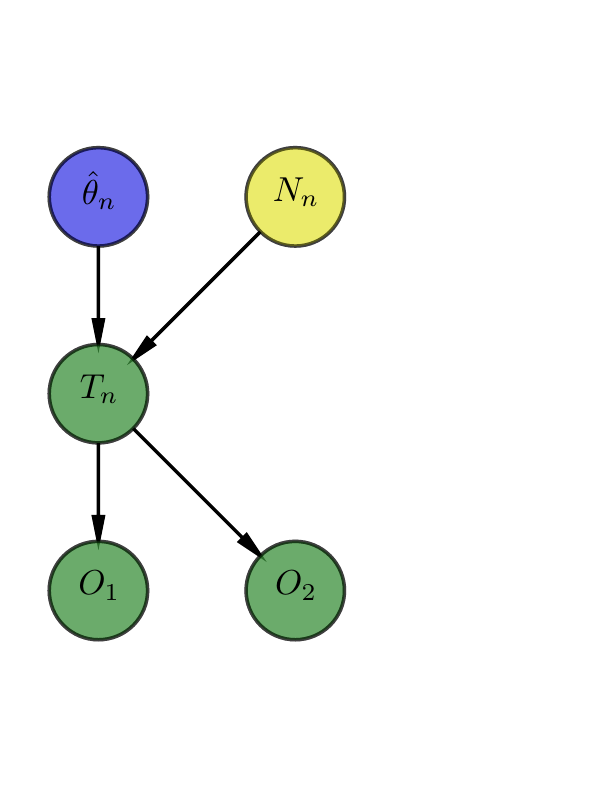}  \\
c)    \includegraphics[width=0.4\textwidth]{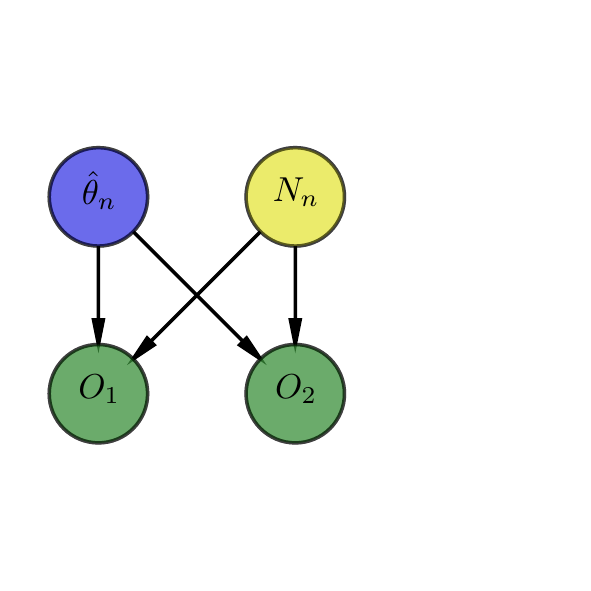} 
    \caption{
a) Application of the selective CLT replacing data $T_n$ with
decomposition \eqref{eq:decomposition}.
\eqref{eq:p0:neyman:carve} has the same dependency graph.
b) The selective sampler applied to a).
c) As $T_n$ is reconstructed directly from $(\hat{\theta}_n, N_n)$ the node $T_n$ can effectively be removed, bearing in mind that it is reconstructed in silico.
}
    \label{fig:CLT}
\end{figure}

If $(T_n, \hat{\theta}_n)$ were jointly normal, then we can decompose
$T_n$ linearly as
\begin{equation} \label{eq:decomposition}
T_n \approx \gamma(F) \cdot \hat{\theta}_n +  N_n
\end{equation}
where $N_n$ ($N$ for nuisance) is asymptotically independent of $\hat{\theta}_n$
$$
N_n = T_n - \frac{\text{Cov}_F(T_n, \hat{\theta}_n)}{\text{Var}_F(\hat{\theta}_n)} \cdot \hat{\theta}_n = T_n - \gamma(F) \cdot \hat{\theta}_n.
$$
When $\theta(F)$ is multivariate, a similar decomposition holds with the reciprocal above replaced by matrix inverse. We note that as long as
$\theta(F)$ is low-dimensional, we will often have reasonably good estimators 
(i.e.~consistent in an appropriate sense) 
of $\text{Var}_F(\hat{\theta}_n)$ under ${\cal M}_n$, for example using
a pairs bootstrap when $T_n$ is formed from IID samples. 
Under certain conditions (c.f.~\cite{randomized_response}) such consistent estimators remain consistent under selective models ${\cal M}_n^*$. 

How does the CLT along ${\cal M}_n$ transfer over to our inferactive setting?
This is illustrated in Figure \ref{fig:CLT}.
Data $T=T_n$ can be replaced by $(\hat{\theta}_n, N_n)$ in the 
dependency graph using the decomposition \eqref{eq:decomposition}
to reconstruct $T_n$ from $(\hat{\theta}_n, N_n)$. The rest of the dependency graph remains the same. This results in a new DAG-DAG given by 
Figure \ref{fig:CLT} a). As in the earlier case, the selective sampler
can be used to reparametrize this DAG-DAG yielding Figure \ref{fig:CLT} b).
Further, as $T_n$ is reconstructed directly from $(\hat{\theta}_n, N_n)$ it can be removed in the DAG-DAG yielding Figure \ref{fig:CLT} c). 

In this final
DAG-DAG, $U$ uses the Gaussian ``plug-in'' estimator from the original
CLT along ${\cal M}_n$. More details, as well as a bootstrap version for the law of $\hat{\theta}_n$ can be found in \cite{randomized_response,selective_bootstrap}. 

We have described a scenario here where $\hat{\theta}_n$ and $N_n$
are jointly Gaussian, though other schemes are also possible. For instance, 
if ${\cal M}_n$ is a sequence of parametric models then $N_n$ can be taken
to be consistent estimates of appropriate nuisance parameters with the
decomposition \eqref{eq:decomposition} replaced by an appropriate kernel.

\subsection{Statistical representation}
\label{sec:repr}

Above, in Section \ref{sec:CLT}, we described 
how existing CLTs can be combined with a DAG-DAG in order to produce
a simpler DAG-DAG in which the statistical model is Gaussian.
Our example only involved one data node $T_n$ with each query being expressed
somehow as a function of $T_n$. In practice, each query may be expressible
as a function of a separate function of $T_n$. For instance, when $T_n=(X_n, Y_n)$ a design matrix and response variable, we see that using the
LASSO to select variables the query can be expressed as a function of $X_n^T Y_n$. 

On the other hand, in the context of our HIV data, suppose $U$ created a boxplot of $Y_n$ against a set of mutations. The vertical demarcations
of the barplot can be expressed
in terms of the median and quantiles of the fold change, stratified across each of these mutations. Such quantiles will typically be subject to a CLT under 
reasonable conditions.
However, the outliers in the boxplot
cannot generally be expressed in terms of quantities that satisfy
a CLT as they are individual sample points.

Each query $i$ then, can be expressed in terms of some random
variable $Z_i$ measurable with respect to the data nodes in the DAG-DAG, 
as well as the randomization that $U$ added to the query. 
We call the function $Z_i$ the {\em statistical representation}
of that query. In silico, the result of a query has many
possible representations that a user might view: a set of non-zero coefficients
may be reported as a list or visualized in an image. This notion
is made explicit in some computing systems, such as the {\tt jupyter} model of 
computation \citep{jupyter}. Our notion of a statistical representation
associates a random variable to the query in such a way that we might
sample from its distribution under our chosen statistical model. 
Many statistical representations are subject to CLT under reasonable (non-selective) statistical models. When the joint distribution
of such representations are asymptotically Gaussian, then there is a clear
way to carry out this program by applying the selective CLT to each query
and conditioning on the appropriate nuisance parameters for each node. Alternatively, one can use the selective bootstrap \citep{selective_bootstrap} to implicitly
carry out this decomposition.

Many model selection queries (LASSO, forward stepwise, marginal screening, and
combinations thereof) have asymptotically Gaussian statistical representations.
An important class of queries, particularly to practicing data scientists, would be graphical queries. Our boxplot example (minus the outlier ticks) is one in which much of the statistical representation is subject to a CLT. 

Scatterplots
are clearly difficult to handle as any scatterplot with the response as one of the variables ``depends'' on all of $Y$ and hence one may have to condition
on the entire response, leaving essentially no information for inference. Randomly perturbing the response is certainly possible. Alternatively, one might try releasing {\em pseudo scatterplots} which fit a simple scatterplot smoother model 
to the data and then report a sample from some parametric model to give $U$ a random sketch of the scatterplot. It seems an interesting question to determine
what features of a scatterplot are actually used by data scientists (and what their statistical representation might be), as well as whether or not
such features actually help the data scientist with their task of discovering
interesting structures in and models for data. We leave such questions for further research.

\subsection{Selection adjusted Bayesian methods}

After observing the results of queries $Q^U_1$ and $Q^U_2$, our data scientist $U$ may have access to scientific literature that might enable him to
choose a parametric model for $T$ along with a prior for its parameters. 
How should the data scientist adjust inference appropriately?

For this, we follow the selection adjusted Bayesian methods of \cite{yekutieli} which were extended to regression models (with an approximate posterior) in \cite{selective_bayesian}. Formally, such methods use the same dependency graph and observed
queries as the frequentist methods though these dependency
graphs are prepended with the appropriate parameters (and possibly hyperparameters). The difference is how a Bayesian data scientist arrives at a posterior. In the above works, the traditional likelihood
that $U$ might use before observing the results of queries
$Q^U_1$ and $Q^U_2$ is replaced with a likelihood truncated to the selection event. Of course, this is nothing more than the selective likelihood
\eqref{eq:likelihood:ratio:1}. Hence, given prior $\pi(\theta)$
the selection adjusted Bayesian methods have posterior
\begin{equation*} \label{eq:bayesian}
	\pi(\theta\:|\:T^i) \:\propto\: \frac{W_{(G_i,Q_i)}(T^i)}{E\left[W_{G_i,Q_i}(T^i)\:|\:\theta\right]} \cdot \pi(\theta).
\end{equation*}

\begin{figure}[h!]
  \centering
a)    \includegraphics[width=0.4\textwidth]{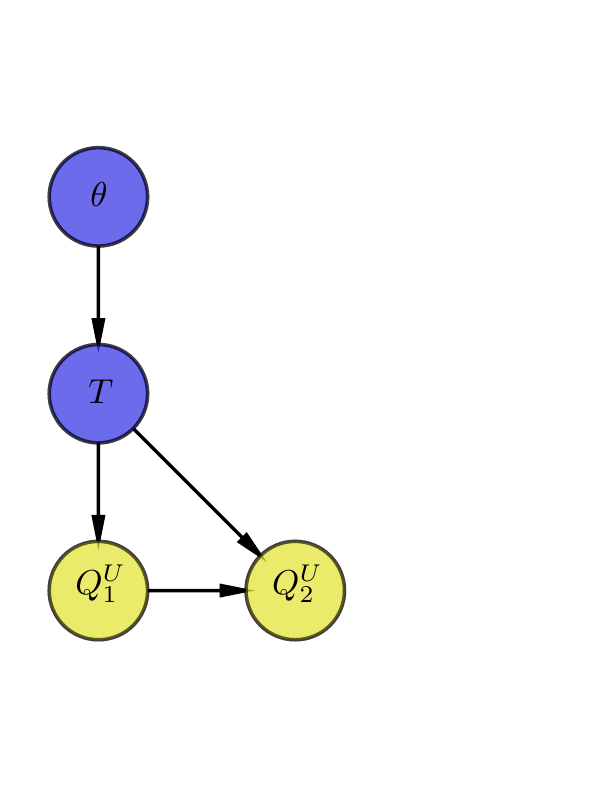}  
b)    \includegraphics[width=0.4\textwidth]{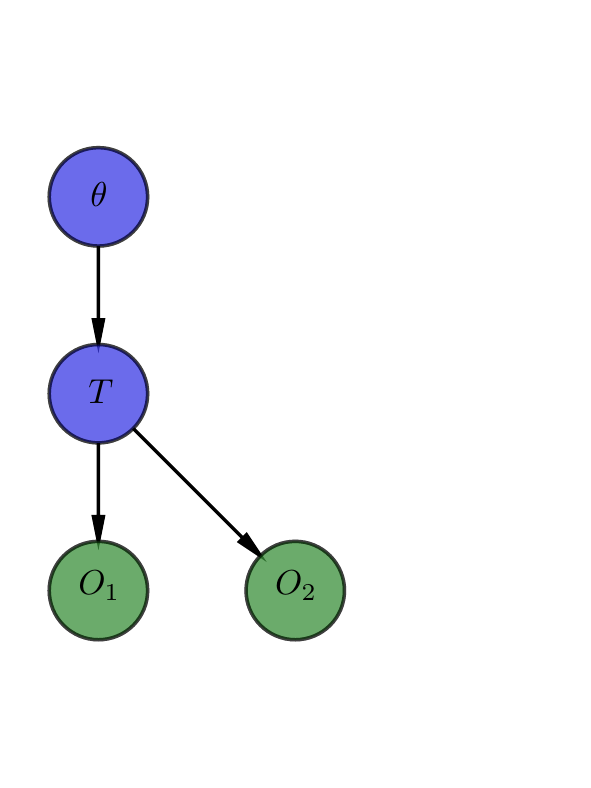}  \\
c)    \includegraphics[width=0.4\textwidth]{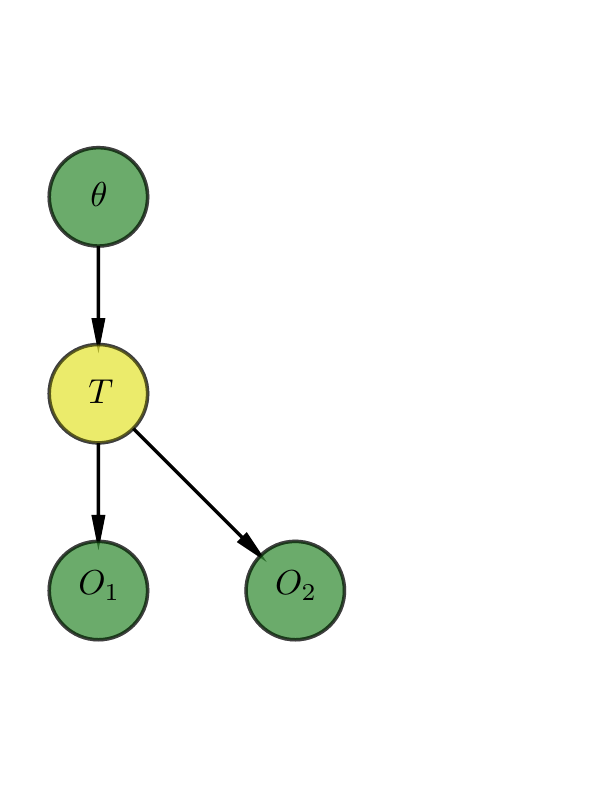} 
    \caption{
a) Selection adjusted Bayesian methods use the same
dependency graph and observed queries with prepended parameters $\theta$.
b) The selective sampler applied to a).
c) For inference, the selection adjusted Bayesian methods condition on the observed data $T$ and marginalizes over everything else, returning the
appropriate posterior.
}
    \label{fig:bayesian}
\end{figure}

The difficulty in using such methods is that the prior depends unavoidably
on $E\left[W_{G_i,Q_i}(T^i)\:|\:\theta\right]$. Approximations based on Chernoff-type bounds
for such probabilities were proposed in \cite{selective_bayesian} to yield tractable
Bayesian inference. The computational cost is still non-negligible: each
step in an MCMC algorithm requires an optimization to approximate this
normalizing constant.
Figure \ref{fig:bayesian} illustrates the construction of the selection adjusted Bayesian posterior.

\section{3TC data revisited} \label{sec:3TC}

In this section, we return to the analysis of the 3TC data and construct
the appropriate reference distribution based on $U$'s data analysis. Let us start with the one-dimensional setting and move on later to scenarios with covariates. We are given data $y=(Y_1,\ldots, Y_n)\overset{i.i.d}{\sim} F_0$ from some distribution $F_0$ with variance 1.
The goal is to do valid inference for $\mu=\mathbb{E}_{Y\sim F_0}[Y]$ with $y$ following the randomized selection event
\begin{equation*}
	\sqrt{n}\bar{y} +\omega >\tau,
\end{equation*}
where $\omega\sim G$ is a random variable in $\mathbb{R}$ independent of $y$. Distribution $G$, with its density denoted as $g$, and threshold $\tau$ are pre-specified and known. The above selection event corresponds to a randomized $z$-test. 
Without selection, we could use the asymptotic normality $\sqrt{n}(\bar{y}-\mu)\overset{d}{\rightarrow}\mathcal{N}(0,1)$ to do inference for the mean. After the randomized selection above, we use the following test statistic
\begin{equation} \label{eq:no:covariates:randomized:pivot}
	\mathcal{P}^R(\sqrt{n}\bar{y};\sqrt{n}\mu,\tau) = (\mathbb{P}_{Z\sim\mathcal{N}(\sqrt{n}\mu,1)}\times\mathbb{P}_{\omega\sim G})\left\{Z<\sqrt{n}\bar{y}\:|\:Z+\omega>\tau\right\}.
\end{equation}
This test statistic is a pivotal quantity, i.e.~it converges to a $\textnormal{Unif}[0,1]$ random variable with respect to the selective distribution of the data, hence it can be used to do valid selective inference. Furthermore, the selective confidence interval for $\mu$ can be constructed by inverting the test $\mathcal{P}^R(\cdot)$. 
The computational and theoretical aspects of the pivot \eqref{eq:no:covariates:randomized:pivot}, including construction of the bootstrap version of the pivot, are given in the appendix. Same techniques apply to regression with covariates as well.\footnote{Inference for the mean in this problem after nonrandomized selection event is also described in detail in the appendix.}

Now we follow the data analysis scenario in the introduction (Section \ref{sec:introduction:example}). Our data $(X, y)$ and the working model is the pairs model that $y_i = \mu (X_i) + \epsilon_i, \ \epsilon_i \sim \mathrm{N} (0, \sigma^2)$ for $i = 1, \ldots, n$ and the conditional mean of $y$ depends on $X$ via the linear predictor $\beta^T X$ where $\beta = (\beta_1, \ldots, \beta_p) \in \mathbb{R}^p$. The goal is to pick the most significant $\beta_i$'s and make valid inference on $\beta$.

Our selective sampler actually enables the selective inference for any convex statistical learning program but suppose the data scientist decides to do a randomized marginal screening thresholded at $c = z_{1-\alpha/2}$ for some nominal $p$-value threshold $\alpha$ on the HIV dataset. Specifically, it takes the $T$ statistics with randomization $\omega_1\in\mathbb{R}^p$ to compute
$$
\tilde{T}_j = T_j + \omega_{1, j} = \frac{X_j^T y}{\hat{\sigma_j} ||X_j||_2} + \omega_{1, j}
$$
for each of the $p$ centered variables $X_j$ and selects the model $E_1 = \{j: |\tilde{T}_j | > c\}$. At this point the data analyst simply reports $E_1$. After we introduce the selective sampler in Section \ref{sec:sampler:marginal:screening}, we will derive valid $p$-values for $E_1$ and also for further multiple queries in Section \ref{sec:sampler:LASSO}.

\subsection{Selective sampler after marginal screening} \label{sec:sampler:marginal:screening}

Following Section \ref{sec: introduction}, the data analyst makes the first query $Q_1^U (t)$ of thresholded marginal screening with randomization which can be expressed in convex optimization form
$$
\minimize_{\eta: \|\eta\|_{\infty} < c}\frac{1}{2}\|T + \omega_1 - \eta\|^2_2, \qquad \omega_1 \sim g_1.
$$
The reconstruction map becomes
\begin{equation*}
	{\phi_1}(T, \eta_{-E_1}, o_{E_1}) = \begin{pmatrix}
		c \cdot s_{E_1} \\ \eta_{-E_1}
	\end{pmatrix} -T+\begin{pmatrix}
		s_{E_1} \cdot o_{E_1} \\ 0
	\end{pmatrix}
\end{equation*}
with the selection event $\mathcal{B}_1$ conditioning on the set achieving the threshold $c$ and their signs to be $(E_1, s_{E_1})$:
\begin{equation*} \label{eq:ms:support}
	\mathcal{B}_1 = \left\{(T,\eta,o_{E_1}): \eta_{E_1} = c\cdot s_{E_1}, \|\eta_{-E_1}\|_{\infty} < c, o_{E_1} \geq 0\right\}.
\end{equation*}
We thus sample $(T,\eta_{-E_1}, o_{E_1})$ from a selective density proportional to  
\begin{equation*} \label{eq:ms:density}
	f(T)\cdot g_1\left(\begin{pmatrix}
		c \cdot s_{E_1} \\ \eta_{-E_1} \end{pmatrix} -T+\begin{pmatrix} s_{E_1} \cdot o_{E_1} \\ 0 \end{pmatrix} \right)
\end{equation*}
and supported on $\mathcal{B}_1$, where $f$ is the unselective law of $T$.

Setting $c = 2.5$ selects 20 mutations and we plot the selective confidence intervals with the observed values for these coefficients (Figure \ref{fig:one:step:randomized:MS}).

\begin{figure}[h!] 
	\begin{center}
  \makebox[\textwidth]{\includegraphics[width=\paperwidth]{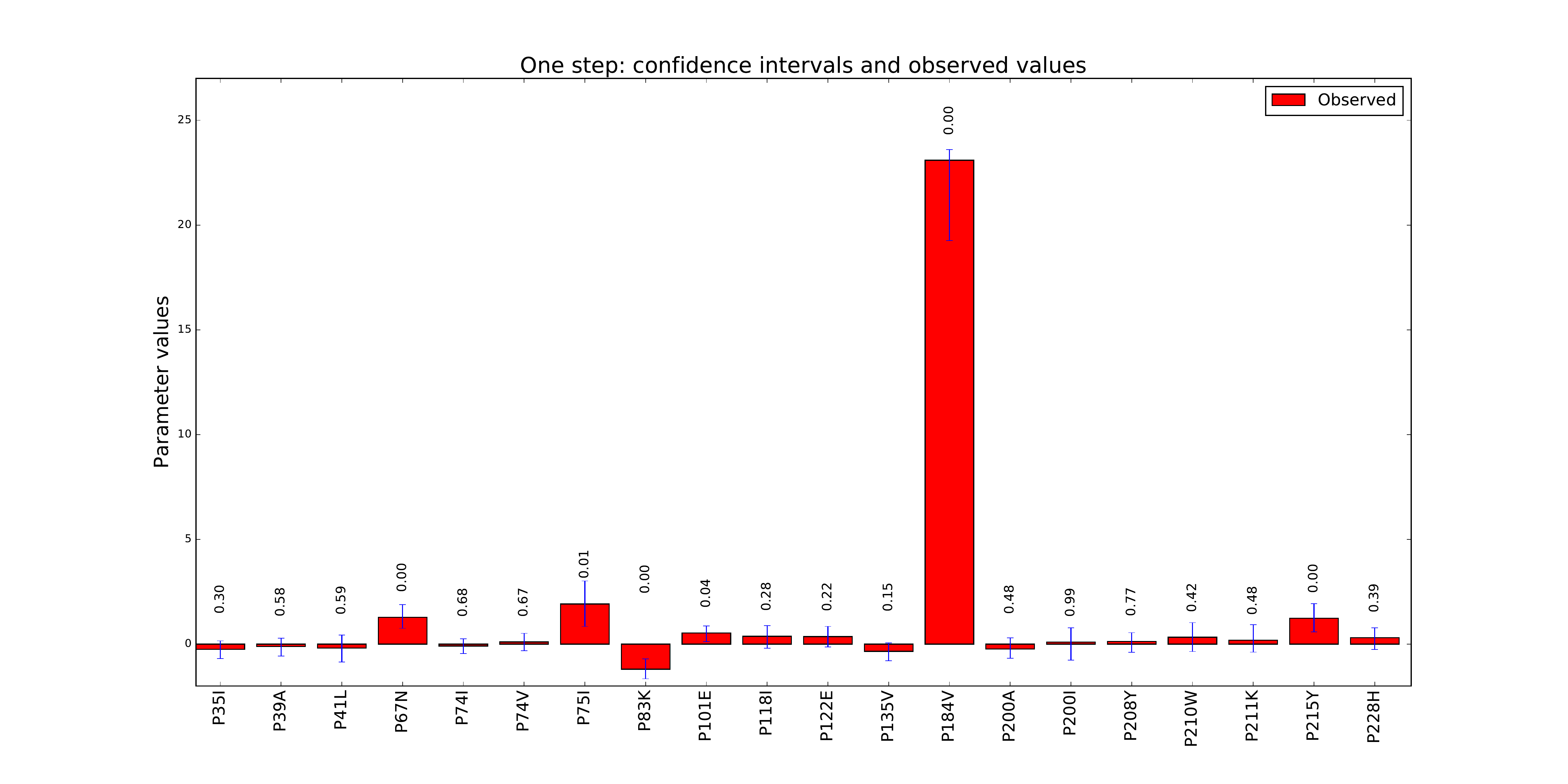}}
\end{center}	
\caption{The plot shows the observed values (red), selective confidence intervals (blue) and selective $p$-values after one view / query on the data.}
\label{fig:one:step:randomized:MS}
\end{figure}

\subsection{Selective sampler after LASSO} \label{sec:sampler:LASSO}

Now the data analyst further fits a LASSO to examine the interactions of selected variables from marginal screening, as proposed by the two-step procedures in \cite{exact_screening}. Specifically, with the selected variables set $E_1$ from the first query $Q_1^U (t)$ we define the (standardized) augmented design matrix $X^{\text{aug}}$ as $E_1$ joined by all the nonzero interactions of them of size $\dbinom{|E_1|}{2}$. Then our second query $Q_2^U (t)$ runs LASSO on $(X^{\text{aug}}, y)$ and selects the support of $\hat{\beta}^{\text{aug}}$ as $E_2 = \{i: \hat{\beta}^{\text{aug}}_i \neq 0\}$. The data analyst could thus compute the selective $p$-value $p_2$ after two queries $(Q_1^U (t), Q_2^U (t))$ of the data.

$Q_2^U (t)$ gives the second optimization problem as
\begin{equation*}  \label{eq:lasso:randomized:program}
	\minimize_{\beta^{\text{aug}} \in \real^p} \frac{1}{2} \|y - X^{\text{aug}}\beta^{\text{aug}}\|^2_2 + \frac{\epsilon}{2} \|\beta^{\text{aug}}\|^2_2 - \omega_2^T \beta^{\text{aug}} + \lambda \|\beta^{\text{aug}}\|_1, \qquad \omega_2 \sim g_2.
\end{equation*}
The reconstruction map for $\omega_2$ is given by
\begin{equation*} 
	{\phi_2}(X^{\text{aug}}, y, \beta^{\text{aug}}_{E_2}, u_{-E_2}) = \epsilon \begin{pmatrix} \beta^{\text{aug}}_{E_2} \\0 \end{pmatrix}-X^{\text{aug}, T} (y-X^{\text{aug}}_{E_2} \beta^{\text{aug}}_{E_2})+\lambda\begin{pmatrix}  s_{E_2} \\ u_{-E_2} \end{pmatrix},
\end{equation*}
with the selection event $\mathcal{B}_2$ conditioning on the active set and their signs to be $(E_2, s_{E_2})$:
\begin{equation*} \label{eq:lasso:randomX:support}
	\mathcal{B}_2 = \left\{(X^{\text{aug}},y, \beta^{\text{aug}}_{E_2},u_{-E_2}): \text{diag}(s_{E_2})\beta^{\text{aug}}_{E_2}> 0, \|u_{-E_2}\|_{\infty} \leq 1 \right\}
\end{equation*}
and the data-dependent Jacobian term $J\psi_2 (X^{\text{aug}},y, \beta^{\text{aug}}_{E_2},u_{-E_2})= \left| \det (X^{\text{aug}, T}_E X^{\text{aug}}_E + \epsilon I) \right|.$
Therefore after $K=2$ steps, the selective density is proportional to 
\begin{equation*}
f_E (X^{\text{aug}}, y) \cdot g_1\left({\phi_1}(T(X, y), \eta_{-E_1}, o_{E_1}) \right) \cdot g_2\left({\phi_2}(X^{\text{aug}}, y, \beta^{\text{aug}}_{E_2}, u_{-E_2}) \right) \cdot J\psi_2 (X^{\text{aug}},y, \beta^{\text{aug}}_{E_2},u_{-E_2})
\end{equation*}
supported on $\left\{ \mathcal{B}_1 (T(X, y), \eta_{-E_1}, o_{E_1}) \cap \mathcal{B}_2 (X^{\text{aug}},y, \beta^{\text{aug}}_{E_2},u_{-E_2})\right\}$. Here $f_E$ denotes the common choice of the reference distribution $f$ where the data analyst sets the support $E = E_1 \cap E_2$. 

Setting $\lambda$ as the theoretical value from \cite{negahban_unified_2012} selects another 12 interaction terms and we plot the selective confidence intervals with the observed values for the coefficients.

\begin{figure}[h!] 
	\begin{center}
  \makebox[\textwidth]{\includegraphics[width=\paperwidth]{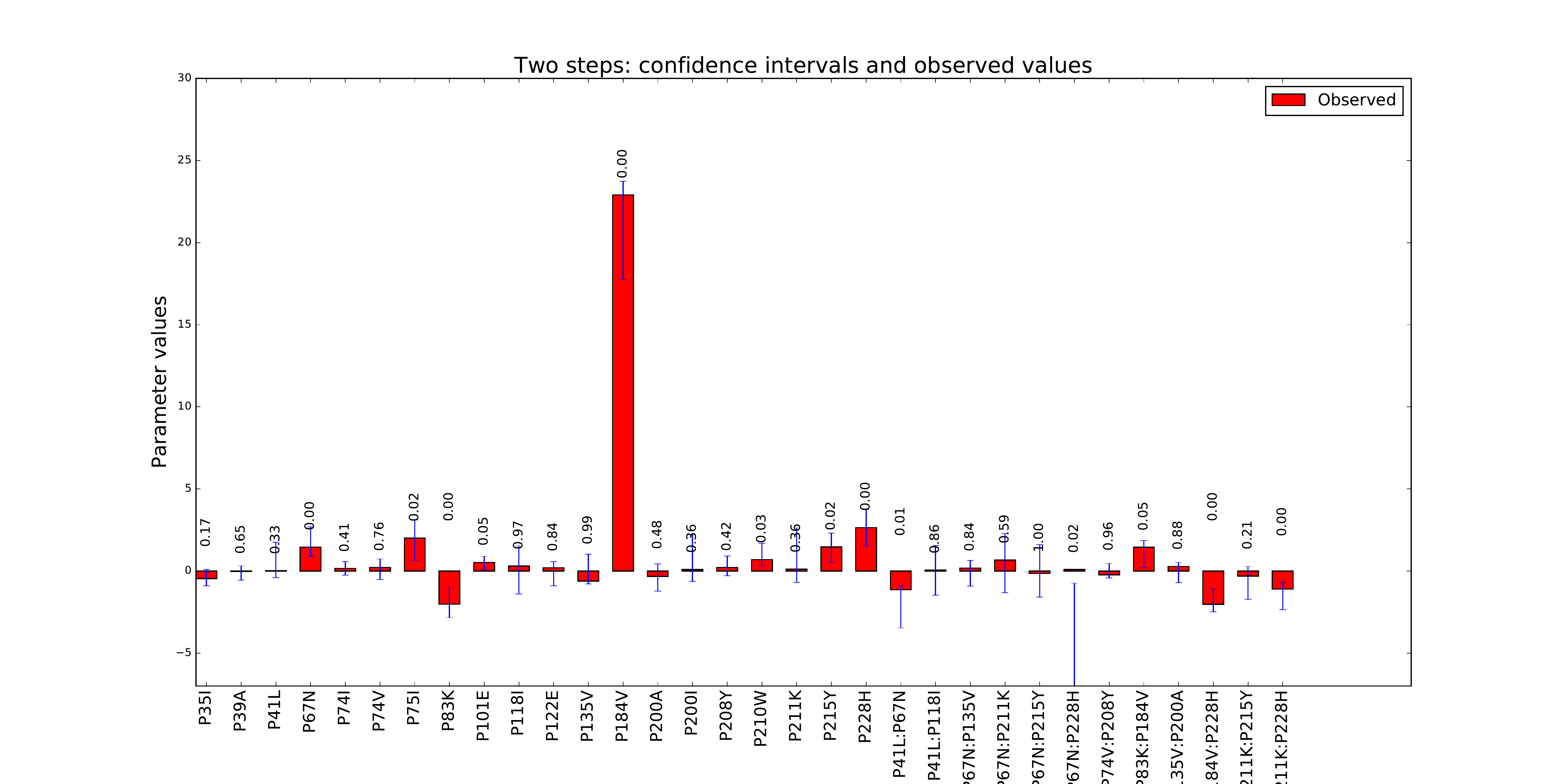}}
\end{center}	
\caption{The plot shows the observed values (red), selective confidence intervals (blue) and selective $p$-values after two views / queries on the data.}
\label{fig:one:step:randomized:MS}
\end{figure}

\section{Discussion}

We have provided a general framework for doing valid inference after making queries on the data and making decisions on what inferential results to report based on the observed outcomes of the queries; we call this ``inferactive data analysis.'' We also introduce a dependency graph, DAG-DAG, to describe the relationship between data and queries. It consists of data and query nodes and can be updated after a data analyst makes additional queries. DAG-DAG becomes very useful when we do selective inference for the chosen parameters, where the conditional density we sample from can be directly obtained from the corresponding graph.

To illustrate how statistical tools are applied in practice, we present various examples in this paper. Among these tools, adding randomization to selection algorithms makes the inferential procedures more powerful and computationally easier. This idea extends to a wide range of popular algorithms, including LASSO, group LASSO, marginal screening, forward-stepwise and their combination into multiple views/queries on the data. For more complicated regression problems, we describe the selection event in terms of a data vector that is asymptotically multivariate Gaussian pre-selection. Theoretically, we can show that the selective CLT and the implied linear decomposition guarantee valid selective inference for the chosen parameters. As we do not use the whole dataset to describe the selection event, we reduce the computational cost.

We illustrate the ``inferactive" procedures through a real HIV dataset. We made two queries on it and presented valid selective p-values and confidence intervals. All the implementations are online.\footnote{\url{https://github.com/jonathan-taylor/selective-inference}.}

\bibliographystyle{agsm}
\bibliography{selective, url=false}

\clearpage
\appendix

\section{Simple problem} \label{sec:simple:problem}

We present an overview of some important results in selective inference literature from the perspective of the ``simple problem.''
In this problem, we are given data $Y_1,\ldots, Y_n\sim F_n$ with mean $\mathbb{E}_{Y\sim F_n}[Y]=\mu(F_n)=\mu_n$ and variance $\mathbb{E}_{Y\sim F_n}[(Y-\mu_n)^2]=\sigma^2 (F_n)=\sigma_n^2$. We assume $\sigma_n^2=1$ in this section but this can be relaxed with the consistent estimate of the variance. Note that $F_n$ and $\mu_n$ can change with $n$. The goal is to do inference for the mean $\mu_n$ after we selected the data vectors $y=(Y_1,\ldots, Y_n)$  for which either nonrandomized or randomized version of the test-statistic $\sqrt{n}\bar{y}$ is greater than a pre-specified threshold. We present the definitions of selective pivots in this case and how they are computed. We further show the constructed pivots are valid post-selection, i.e.~using them for inference for $\mu_n$ ensures the (asymptotic) control of selective type I error.
 These results translate to regression problems with more complicated affine selection events.

Let us define the following notations.
\begin{itemize}[leftmargin=*]
	\item $F_n^n$ denotes the distribution of $n$ IID samples from $F_n$.
	\item $\hat{F}_n$ denotes the empirical distribution of $y$ and as above $\hat{F}_n^n$ denotes the distribution of $n$ IID samples from $\hat{F}_n$.
	\item $F_n^*$ denotes the joint distribution of the data vector $y$ after selection, either randomized or nonrandomized.
\end{itemize}

\subsection{Simple example: nonrandomized selection event}

\subsubsection{Nonrandomized pivots: definitions and computations}

The nonrandomized selection event is of the form
\begin{equation*}
	\sqrt{n}\bar{y}>\tau_n,
\end{equation*}
for a given threshold $\tau_n$. We show how to do inference for $\mu_n$ after our data $y$, originally generated from $F_n^n$, has been selected to satisfy the above inequality.

Before defining the selective pivots for $\mu_n$, let us define two conditional CDFs. 
\begin{itemize}[leftmargin=*]
	\item If the test statistic $\sqrt{n}\bar{y}$ was normal, i.e.~from $\mathcal{N}(\sqrt{n}\mu_n,1)$, then the conditional CDF of $\sqrt{n}\bar{y}$ would be
	\begin{equation*}
		F^N_n(t;\sqrt{n}\mu_n,\tau_n) = \mathbb{P}_{Z\sim\mathcal{N}(\sqrt{n}\mu_n,1)}\left\{Z<t\:|\:Z>\tau_n\right\}
	\end{equation*}
	for all $t\in\mathbb{R}$. 
	\item If the data was coming from its empirical distribution, then the conditional distribution of $\sqrt{n}\bar{y}$ would be
		\begin{equation*}
			F^B_n(t;\sqrt{n}\mu_n,\tau_n)=\mathbb{P}_{y^{\star}\sim\hat{F}_n^n}\left\{ \sqrt{n}(\bar{y}^{\star}-\bar{y})+\sqrt{n}\mu_n <t \:|\:\sqrt{n}(\bar{y}^{\star}-\bar{y})+\sqrt{n}\mu_n>\tau_n\right\}
		\end{equation*}
	for all $t\in\mathbb{R}$, where the RHS is with respect to the bootstrap sample $y^{\star}=(Y_1^{\star},\ldots, Y_n^{\star})\sim\hat{F}_n^n$.
\end{itemize}


The selective pivots for this problem are defined and computed as follows.
\begin{enumerate}[leftmargin=*]
	\item \textit{Plugin Gaussian pivot after nonrandomized selection or TG pivot} is defined as
		\begin{equation*}
			\mathcal{P}^N(\sqrt{n}\bar{y};\sqrt{n}\mu_n,\tau_n)
			=F_n^N(\sqrt{n}\bar{y};\sqrt{n}\mu_n,\tau_n) 
			=\frac{\Phi(\sqrt{n}(\bar{y}-\mu_n))-\Phi(\tau_n-\sqrt{n}\mu_n)}{1-\Phi(\tau_n-\sqrt{n}\mu)},
		\end{equation*}
		where $\Phi$ denotes the CDF of a standard normal distribution.
		This pivot is a truncated Gaussian (TG) test statistic introduced in \cite{exact_lasso}. 
		 Since this pivot has an explicit form, computing it is easy. 

	
	\item \textit{Bootstrap pivot after nonrandomized selection} is defined as
		\begin{equation*}
		\begin{aligned}
			&\mathcal{P}^B(\sqrt{n}\bar{y};\sqrt{n}\mu_n,\tau_n)
			=F^B_n(\sqrt{n}\bar{y};\sqrt{n}\mu_n,\tau_n) \\
			&\qquad=\frac{\hat{F}_n^n\left\{\sqrt{n}(\bar{y}^{\star}-\bar{y})<\sqrt{n}(\bar{y}-\mu_n)\right\}-\hat{F}_n^n\left\{\sqrt{n}(\bar{y}^{\star}-\bar{y})<\tau_n-\sqrt{n}\mu_n\right\}}{1-\hat{F}^n_n\left\{\sqrt{n}(\bar{y}^{\star}-\bar{y})<\tau_n-\sqrt{n}\mu\right\}}.
		\end{aligned}
		\end{equation*}
		This pivot is introduced in \cite{cmu_bootstrap}.  To compute the bootstrap pivot above we sample data with replacement to get the bootstrap samples $(Y_1^{\star b},\ldots, Y_n^{\star b})$ for $b=1,\ldots,B$, where the size $B$ is large. We then compute the pivot as
		\begin{equation*}
			\frac{\frac{1}{B}\sum_{b=1}^B\mathbb{I}_{\{\sqrt{n}(\bar{y}^{\star b}-\bar{y})<\sqrt{n}(\bar{y}-\mu_n)\}}
			-\frac{1}{B}\sum_{b=1}^B\mathbb{I}_{\{\sqrt{n}(\bar{y}^{\star b}-\bar{y})<\tau_n-\sqrt{n}\mu_n\}}}{1-\frac{1}{B}\sum_{b=1}^B\mathbb{I}_{\{\sqrt{n}(\bar{y}^{\star b}-\bar{y})<\tau_n-\sqrt{n}\mu_n\}}}.
		\end{equation*}
\end{enumerate}

Note that we can easily invert either of these pivots to get a confidence interval for $\mu_n$ by evaluate the pivot over a grid of parameter values or use some of the root-finding functions.


\subsubsection{Asymptotics of the nonrandomized pivots and honest coverage post-selection}

We now state the results showing that under the post-selective distribution the pivots defined above have asymptotically Unif$[0,1]$ distribution. We denote as $F_n^*$ the post-selective distribution of the data $y\sim F^n_n$ conditional on selection. 
Given $\{F_n\in\mathcal{F}_n:n\geq 1\}$ as a sequence of distribution classes $\{\mathcal{F}_n:n\geq 1\}$ to prove convergence results we state the following assumptions that we refer to separately.
\begin{enumerate}[label=(\Alph*),leftmargin=3em]
	\item \label{item:param:conv} Convergence of the parameters: for every sequence of distributions $\{F_n: F_n\in\mathcal{F}_n, n\geq 1\}$ there exists $\mu$ such that $\sqrt{n}\mu(F_n)\rightarrow\mu$ as $n\rightarrow\infty$.
	\item \label{item:second:moment} Second moment assumption: for all $F_n \in\mathcal{F}_n$ and all $n\geq 1$, $\mathbb{E}_{X\sim F_n}[X^2]=1$.
	\item \label{item:third:moment} Third moment assumption: for all $F_n \in\mathcal{F}_n$ and all $n\geq 1$, $\mathbb{E}_{X\sim F_n}[X^3]<\infty$.
\end{enumerate}

We will also need the following assumption about the selection regions.
\begin{enumerate}[label=(\Alph*), leftmargin=3em]
\setcounter{enumi}{3}
	\item \label{item:tau:conv} Convergence of the selection region: $\tau_n\rightarrow\tau$ as $n\rightarrow\infty$.
\end{enumerate}

\begin{lemma}[Nonrandomized selection: asymptotics of the plugin Gaussian pivot] \label{lemma:non:randomized:plugin:gaussian}
	Assuming $\{\mathcal{F}_n:n\geq 1\}$ is a sequence of distributions for which \ref{item:param:conv} and \ref{item:second:moment} hold. Further assuming \ref{item:tau:conv} holds, we have
	\begin{equation*}
		\underset{n\rightarrow\infty}{\lim}\:\underset{F_n\in\mathcal{F}_n}{\sup}\:\underset{t\in[0,1]}{\sup}\left|F_n^*\{\mathcal{P}^N(\sqrt{n}\bar{y};\sqrt{n}\mu_n,\tau_n)\leq t\} -t\right|=0.\footnote{Note that $\mathbb{F}_n^*$ only depends on $\mathbb{F}_n$ and $\tau_n$ so it suffices to write the supremum over $\mathcal{F}_n$.}
	\end{equation*}
	Furthermore, if we define $CI(\sqrt{n}\bar{y}; \tau_n,\alpha)$ to be the set of $\mu_n$ such that $\alpha/2\leq \mathcal{P}(\sqrt{n}\bar{y};\sqrt{n}\mu_n,\tau_n)\leq 1-\alpha/2$, then $CI^N(\sqrt{n}\bar{y};\tau_n,\alpha)$ is asymptotically honest for $\mu_n$, i.e.
	\begin{equation*}
		\underset{n\rightarrow\infty}{\lim}\:\underset{F_n\in\mathcal{F}_n}{\sup}\:\underset{\alpha\in[0,1]}{\sup}\left|F_n^*\{\mu_n\in CI^N(\sqrt{n}\bar{y};\tau_n,\alpha)\}-(1-\alpha)\right| =0.
	\end{equation*}
\end{lemma}

\begin{lemma}[Nonrandomized selection: asymptotics of the bootstrap pivot] \label{lemma:non:randomized:boot}
	Assuming $\{\mathcal{F}_n:n\geq 1\}$ is a sequence of distributions for which  \ref{item:param:conv}, \ref{item:second:moment} and \ref{item:third:moment} hold. Further assuming  \ref{item:tau:conv}  holds, we have
	\begin{equation*}
		\underset{n\rightarrow\infty}{\lim}\:\underset{F_n\in\mathcal{F}_n}{\sup}\:\underset{t\in[0,1]}{\sup}\left|F_n^*\{\mathcal{P}^B(\sqrt{n}\bar{y};\sqrt{n}\mu_n,\tau_n)\leq t\} -t\right|=0.
	\end{equation*}
	Furthermore, if we define $CI^B(\sqrt{n}\bar{y}; \tau_n, \alpha)$ to be the set of $\mu_n$ such that $\alpha/2\leq \mathcal{P}^B(\sqrt{n}\bar{y};\sqrt{n}\mu_n,\tau_n)\leq 1-\alpha/2$, then $CI^B(\sqrt{n}\bar{y};\tau_n,\alpha)$ is asymptotically honest for $\mu_n$, i.e.
	\begin{equation*}
		\underset{n\rightarrow\infty}{\lim}\:\underset{F_n\in\mathcal{F}_n}{\sup}\:\underset{\alpha\in[0,1]}{\sup}\left|F_n^*\{\mu_n\in CI^B(\sqrt{n}\bar{y};\tau_n,\alpha)\}-(1-\alpha)\right| =0.
	\end{equation*}
\end{lemma}

Note that the results above are under the conditional distribution $F_n^*$. \cite{tian2015asymptotics} show that the TG pivot after LASSO is asymptotically uniform with a different set of assumptions. For low-dimensional regression problems, if the selection event is affine in $y$ (or $T$), \cite{markovic2017adaptive} gives a more general version of Lemma \ref{lemma:non:randomized:plugin:gaussian}, proving the result with $F_n^*$. \cite{cmu_bootstrap} proves the result with respect to unconditional distribution $F_n$. \cite{cmu_bootstrap} further proves that bootstrap pivot leads to asymptotically conservative confidence intervals.

\subsubsection{Post-selection consistency result}

Given that $\bar{y}$ is consistent pre-selection (under $F_n$) for $\mu_n$, the following lemma shows that $\bar{y}$ is also consistent post-selection (under $F_n^*$) for $\mu_n$ under some conditions. Note that this results generalizes to any consistent estimators but to keep the notation simpler we take the estimator to be $\bar{y}$ and the parameter to be $\mu_n$.

\begin{lemma}[Nonrandomized selection: post-selection consistency of the mean] 
	We are given a sequence of distribution classes $\{\mathcal{F}_n:n\geq 1\}$ over which
	 $\bar{y}$ is consistent pre-selection, i.e.~for every $\epsilon>0$, $\underset{n\rightarrow\infty}{\lim}\underset{F_n\in\mathcal{F}_n}{\sup}F_n\left\{|\bar{y}-\mu_n|>\epsilon\right\}=0$. Let us further assume that \ref{item:param:conv} and \ref{item:tau:conv} hold. Then, $\bar{y}$ is consistent post-selection as well, i.e.
	 \begin{equation*}
	 	\forall\epsilon, \;\;\;\underset{n\rightarrow\infty}{\lim}\:\underset{F_n\in\mathcal{F}_n}{\sup}F_n^*\left\{|\bar{y}-\mu_n|>\epsilon\right\}=0.
	 \end{equation*}
\end{lemma}


\subsection{Simple example: randomized selection event}

\subsubsection{Randomized pivots: definitions and computations}

As above, our goal is to do inference for $\mu_n$ after the randomized selection event
\begin{equation*}
	\sqrt{n}\bar{y}+\omega>\tau_n,
\end{equation*}
where the randomization $\omega\sim G$ is independent of $y$. The distribution $G$, with density denoted as $g$, is determined by the user.

Before defining the pivots in this case, let us define two CDFs that marginalize over $\omega$.
\begin{itemize}[leftmargin=*]
\item If the data was Gaussian, i.e.~$\sqrt{n}\bar{y}\sim\mathcal{N}(\sqrt{n}\mu_n,1)$, the conditional distribution of the data would be
\begin{equation*}
	F^{N,R}_n(t;\sqrt{n}\mu_n,\tau_n) = (\mathbb{P}_{Z\sim\mathcal{N}(\sqrt{n}\mu_n,1)}\times \mathbb{P}_{\omega\sim G})\left\{Z<t\:|\:Z+\omega>\tau_n\right\}.
\end{equation*}
\item If the data was generated from its empirical distribution, the conditional distribution of the data would be
\begin{equation*}
	F^{B,R}_n(t;\sqrt{n}\mu_n,\tau_n) = (\mathbb{P}_{y^{\star}\sim\hat{F}_n^n}\times \mathbb{P}_{\omega\sim G})\left\{\sqrt{n}(\bar{y}^{\star}-\bar{y})+\sqrt{n}\mu_n<t\:|\:\sqrt{n}(\bar{y}^{\star}-\bar{y})+\sqrt{n}\mu_n+\omega>\tau_n\right\}.
\end{equation*}
\end{itemize}

The pivots in these cases are defined and computed as follows.
\begin{enumerate}[leftmargin=*]
	\item \textit{Plugin Gaussian pivot after randomized selection event} is defined as
	\begin{equation*}
		\mathcal{P}^{N,R}(\sqrt{n}\bar{y};\sqrt{n}\mu_n, \tau_n) =  	F^{N,R}_n(\sqrt{n}\bar{y}; \sqrt{n}\mu_n, \tau_n).
 	\end{equation*} 
	This pivot is introduced in \cite{randomized_response}. 
	It can be computed in following two ways.
	
	\begin{enumerate}[leftmargin=*]
	\item\label{item:plugin:gaussian:sampling} \textit{Sampling after the change of variables.} The post-selection distribution of $(\sqrt{n}\bar{y},\omega)$ pair is proportional to
	\begin{equation*}
		f(\sqrt{n}\bar{y})\cdot g(\omega)\cdot \mathbb{I}_{\{\sqrt{n}\bar{y}+\omega>\tau_n\}},
	\end{equation*}
	where $f$ is the density of $\sqrt{n}\bar{y}$ pre-selection. Based on the pre-selection asymptotic normality $\sqrt{n}(\bar{y}-\mu_n)\overset{d}{\rightarrow}\mathcal{N}(0, 1)$ as $n\rightarrow\infty$, where the LHS is under $y\sim F_n^n$, we can assume $f$ is $\phi_{(\sqrt{n}\mu_n,1)}$, the density of $\mathcal{N}(\sqrt{n}\mu_n,1)$.
	The selective sampler of \cite{selective_sampler} avoids sampling from the density above that has the constraints on a linear combination of data and randomization. Its idea is to do a change of variables $z=\sqrt{n}\bar{y}+\omega$. In this case the post-selection density in terms of $(\sqrt{n}\bar{y},z)$ becomes proportional to
	\begin{equation} \label{eq:selective:density:not-marginalized}
		f(\sqrt{n}\bar{y})\cdot g(z-\sqrt{n}\bar{y})\cdot\mathbb{I}_{\{z>\tau_n\}}.
	\end{equation}
	The variable $z$ is called the optimization variable.
	Note that the density above does not have any restrictions on the data $\sqrt{n}\bar{y}$ but only on $z$. We can sample $(\sqrt{n}\bar{y},z)$ from the density above by moving both $(\sqrt{n}\bar{y},z)$ using MCMC. In practice, we often use projected Langevin to sample.
	In this case, we can further marginalize over the optimization variable $z$ and get the post-selection density of the data as proportional to
	\begin{equation} \label{eq:selective:density:marginalized}
		f(\sqrt{n}\bar{y}) \cdot \int\limits_{z=\tau_n}^{\infty} g(z-\sqrt{n}\bar{y}) dz = f(\sqrt{n}\bar{y})\cdot\left(1-G(\tau_n-\sqrt{n}\bar{y})\right).
	\end{equation}
	For more complicated selection events we sample the data and the optimization variables from a selective density \citep{selective_sampler}.

	\item \label{item:plugin:gaussian:approx} \textit{Approximation of selection probabilities.} In this example, we sample $\sqrt{n}\bar{y}$ from the density in \eqref{eq:selective:density:marginalized} or $(\sqrt{n}\bar{y},z)$ from the density in \eqref{eq:selective:density:not-marginalized}. Another way to get the post-selective distribution of the data is to evaluate the expression from \eqref{eq:selective:density:marginalized}  written as $f(t)\cdot(1-G(\tau_n-t))$ over the grid of $t$ values, where we use $t$ to denoted the argument in the density corresponding to random variable $\sqrt{n}\bar{y}$. Based on the set of these values over a grid, we can approximate the normalizing constant to get the selective density of the data. This approach has been developed for more complicated randomized selection events, including LASSO, $\ell_1$-penalized regression, first step of forward-stepwise and marginal screening \citep{panigrahi2017mcmc}. This is the same as that we weight the pre-selection density $f(t)$ with $1-G(\tau_n-t)=\mathbb{P}_{\omega\sim G}\{t+\omega>\tau_n\}$, which is the selection probability.  For more complicated selection events we do not have these probabilities over a range of $t$ values explicitly but we use approximation of \cite{selective_bayesian} to compute them approximately.
	\end{enumerate}
	
		
	\item \textit{Bootstrap pivot after randomized selection event} is defined as
	\begin{equation*}
		\mathcal{P}^{B,R}(\sqrt{n}\bar{y};\sqrt{n}\mu_n, \tau_n) =  	F^{B,R}_n(\sqrt{n}\bar{y};\sqrt{n}\mu_n, \tau_n).
 	\end{equation*}
	This pivot is introduced in \cite{selective_bootstrap}. In this case there are two ways to compute the bootstrap pivot.
	
	\begin{enumerate}[leftmargin=*]
	\item \label{item:boot:pivot:wild}\textit{Wild bootstrap.} 
	Standard wild bootstrap without selection approximates the distribution of $\sqrt{n}(\bar{y}-\mu_n)$ with the distribution of $\frac{1}{\sqrt{n}}\sum_{i=1}^n(Y_i-\mu_n)\alpha_i$ where $\bm\alpha=(\alpha_1,\ldots,\alpha_n)\overset{IID}{\sim}\mathcal{N}(0,1)$ are bootstrap weights  taken independently of the data.
	After selection, we write the post-selection distribution of $(z,\bm\alpha)$ as proportional to
	 \begin{equation} \label{eq:wild:boot:density}
	 	\left(\prod_{i=1}^n\phi_{(0,1)}(\alpha_i)\right)\cdot g\left(z-\frac{1}{\sqrt{n}}\sum_{i=1}^n(Y_i-\mu_n)\alpha_i+\sqrt{n}\mu_n\right)\cdot\mathbb{I}_{\{z>\tau_n\}},
	 \end{equation}
	where $z$ is the optimization variable defined above. We sample from the density in \eqref{eq:wild:boot:density} using MCMC. This approach has been developed in \cite{selective_bootstrap}.

	\item \label{item:boot:pivot:standard} \textit{Sampling with replacement.} Given $B$ bootstrap samples $(Y_1^{\star b},\ldots,Y_n^{\star b})\overset{IID}{\sim}\hat{F}_n$, for all $b=1,\ldots, B$, the bootstrap pivot is computed as
		\begin{equation*}
		\begin{aligned}
			&\mathcal{P}^{B,R}(\sqrt{n}\bar{y}; \sqrt{n}\mu_n, \tau_n) \\
			&=\frac{\frac{1}{B}\sum_{b=1}^B \mathbb{P}_{\omega\sim G}\left\{\sqrt{n}(\bar{y}^{\star b}-\bar{y})+\sqrt{n}\mu_n+\omega>\tau_n\right\}\mathbb{I}_{\{\sqrt{n}(\bar{y}^{\star b}-\bar{y})\leq \sqrt{n}(\bar{y}-\mu_n)\}}}{\frac{1}{B}\sum_{b=1}^B \mathbb{P}_{\omega\sim G}\left\{\sqrt{n}(\bar{y}^{\star b}-\bar{y})+\sqrt{n}\mu_n+\omega>\tau_n\right\}} \\
			&= \frac{\frac{1}{B}\sum_{b=1}^B\left(1- G(\tau_n-\sqrt{n}(\bar{y}^{\star b}-\bar{y})-\sqrt{n}\mu_n)\right)\mathbb{I}_{\{\sqrt{n}(\bar{y}^{\star b}-\bar{y})\leq\sqrt{n}(\bar{y}-\mu_n)\}}}{1-\frac{1}{B}\sum_{b=1}^B G(\tau_n-\sqrt{n}(\bar{y}^{\star b}-\bar{y})-\sqrt{n}\mu_n)}.
		\end{aligned}
		\end{equation*}
		In words, to compute the bootstrap pivot we weight each of the bootstrap samples $\sqrt{n}(\bar{y}^{\star b}-\bar{y})$ with the probability of the event $\{ \sqrt{n}(\bar{y}^{\star b}-\bar{y})+\sqrt{n}\mu_n+\omega>\tau_n \}$ marginalized over $\omega\sim G$ only. In these examples as well, we weight the bootstrap samples with the probabilities of selection corresponding to those samples where the probability as above is over the randomization. These probabilities are computed approximately using the barrier approximation of \cite{selective_bayesian} and the results will be presented in future work.
		\end{enumerate}

	\end{enumerate}


To invert the pivots from either \eqref{item:plugin:gaussian:sampling} or \eqref{item:boot:pivot:wild} and construct a confidence interval for $\mu_n$ we sample once using some estimated $\mu_n$ value as a reference parameter (usually MLE pre-selection). Then we tilt the samples to get the pivots evaluated at other parameter values. To invert the pivots from either \eqref{item:plugin:gaussian:approx} or \eqref{item:boot:pivot:standard}, we evaluate the pivots separately over a grid of parameter values.

\subsubsection{Asymptotics of the randomized pivots and the honest coverage post-selection}

We now state the results showing that under the post-selective distribution the pivots $\mathcal{P}^R$ and $\mathcal{P}^{R,B}$ defined above have asymptotically Unif$[0,1]$ distribution. We denote as $F_n^{*R}$ the post-selective distribution of the data $y\sim F^n_n$ conditional on randomized selection, marginalizing over the randomization. Also note that the stated asymptotic results about the pivots imply the confidence intervals constructed are uniformly honest over the corresponding classes. We do not state the coverage results separately but the guarantees for them are the same as the ones given in Lemma \ref{lemma:non:randomized:plugin:gaussian} and Lemma \ref{lemma:non:randomized:boot}.

Before stating the lemmas let us introduce two more assumptions. The first one is on the distance between the mean parameter $\mu_n$ and the selection region. The second one assumes the distribution of the randomization has heavier tails than Gaussian.
\begin{enumerate}[label=(\Alph*), leftmargin=3em]
\setcounter{enumi}{4}
	\item \label{item:local:alternatives} Local alternatives assumption: there exists a constant $B$ such that for every sequence $\{F_n:F_n\in\mathcal{F}_n, n\geq 1\}$ we have that the distance between the parameter $\sqrt{n}\mu_n$ and the selection region $[\tau_n,\infty)$ is not too big, i.e.~$\max\{\tau_n-\sqrt{n}\mu(F_n)),0\}<\Delta$, for some constant $\Delta$.
	\item \label{item:lip:randomization} Lipschitz assumption on the randomization: Assume the density $g(x)$, $x\in\mathbb{R}$, is proportional to $\exp(-\widetilde{g}(x))$, $x\in\mathbb{R}$, where $\widetilde{g}$ is a Lipschitz continuous function with smooth derivatives up to the third order.
\end{enumerate}

Depending on whether the randomization is Gaussian or with heavier tails,  we need different assumptions for the asymptotic convergence of the plugin pivots, and thus we state these two results separately. More explicitly, we do not need any extra assumption for Lipschitz randomization (i.e.~we allow for $\mu_n$ to be arbitrarily far from the selection region, implying that the constructed pivot is valid even after rare selection events). On the other hand, we need \ref{item:local:alternatives} for Gaussian randomization.

\begin{lemma}[Randomized selection: asymptotics of the plugin Gaussian pivot] \label{lemma:randomized:plugin:gaussian}
Assume either of the two sets of assumptions hold for the sequence $\{\mathcal{F}_n:n\geq 1\}$, selection region and the randomization distribution $G$:
\begin{enumerate}
	\item \ref{item:second:moment}, \ref{item:local:alternatives} and $G=\mathcal{N}(0,\gamma^2)$ or
	\item \ref{item:second:moment} and \ref{item:lip:randomization}.
\end{enumerate}
Then we have
	\begin{equation*}
		\underset{n\rightarrow\infty}{\lim}\:\underset{F_n\in\mathcal{F}_n}{\sup}\:\underset{t\in[0,1]}{\sup}\left|F_n^*\{\mathcal{P}^{N,R}(\sqrt{n}\bar{y};\sqrt{n}\mu_n,\tau_n)\leq t\} -t\right|=0.
	\end{equation*}
\end{lemma}

\begin{lemma}[Randomized selection: asymptotics of the bootstrap pivot]
\label{lemma:randomized:boot}
	Assuming $\{\mathcal{F}_n:n\geq 1\}$ is such that the assumptions \ref{item:second:moment} and \ref{item:third:moment} above hold. Further assume \ref{item:local:alternatives} and \ref{item:lip:randomization} hold. Then we have 
	\begin{equation*}
		\underset{n\rightarrow\infty}{\lim}\:\underset{F_n\in\mathcal{F}_n}{\sup}\:\underset{t\in[0,1]}{\sup}\left|F_n^*\{\mathcal{P}^{R,B}(\sqrt{n}\bar{y};\sqrt{n}\mu_n,\tau_n)\leq t\} -t\right|=0.
	\end{equation*}
\end{lemma}

A more general version of Lemma \ref{lemma:randomized:plugin:gaussian} for affine selection regions is first derived in \cite{tian2015asymptotics} under both local alternatives  and heavy-tailed randomization assumptions. \cite{selective_bootstrap} extends this result by removing the local alternatives assumption for heavy tailed randomizations and proves the general versions of Lemma \ref{lemma:randomized:plugin:gaussian} and Lemma \ref{lemma:randomized:boot} for affine selection regions.

\subsubsection{Post-selection consistency results}

Considering randomized selection event, we present some results here on post-selection consistency of estimators. By the strong law of large numbers, $\bar{y}$ is consistent for the mean $\mu_n$ pre-selection, i.e.~with respect to the pre-selection distribution of the data $F_n$. We state the results of \cite{randomized_response, selective_bootstrap} showing that the sample mean $\bar{y}$ is consistent for the mean $\mu_n$ post-selection, i.e.~with respect to $F_n^*$. In the case of Gaussian data and Gaussian randomization, we present a result of \cite{selective_bayesian} showing a different estimator, called selective MLE, is consistent post-selection.

Let us first state the results about the post-selection consistency of $\bar{y}$. The case of Gaussian randomization requires local alternatives assumption while the case of Lipschitz randomization does not. These results generalize to any consistent estimators \citep{randomized_response, selective_bootstrap}.

\begin{lemma}[Randomized selection: post-selection consistency of the mean] 
	We are given a sequence of distributions $\{\mathcal{F}_n:n\geq 1\}$ over which
	 $\bar{y}$ is consistent pre-selection, i.e.~for every $\epsilon>0$ we have $\underset{n\rightarrow\infty}{\lim}\underset{F_n\in\mathcal{F}_n}{\sup}F_n\left\{|\bar{y}-\mu_n|>\epsilon\right\}=0$.
	 Further assume either of the two following sets of assumptions holds:
	 \begin{enumerate}
	 	\item $G=\mathcal{N}(0,1)$ and \ref{item:local:alternatives} or
	 	\item \ref{item:lip:randomization}.
	 \end{enumerate}
	Then we have that the mean is consistent post-selection as well, i.e.
	 \begin{equation*}
	 	\forall\epsilon,\;\;\;\underset{n\rightarrow\infty}{\lim}\:\underset{F_n\in\mathcal{F}_n}{\sup}F_n^*\left\{|\bar{y}-\mu_n|>\epsilon\right\}=0.
	 \end{equation*}
\end{lemma}

In the case of Gaussian randomization, local alternatives assumption is necessary for $\bar{y}$ to be consistent for the mean $\mu_n$; a counter example was given in \cite{selective_bayesian} without the assumption. In such scenarios when the local alternative assumption does not hold, \cite{selective_bayesian} proposes another estimator, called selective MLE, and they show it is consistent for the mean.

In order to define  the selective MLE, let us derive the post-selection density of the data assuming the pre-selection distributions of the data and randomization are normal, i.e.~$F_n=\mathcal{N}(\mu_n,1)$ and $G=\mathcal{N}(0,\gamma^2)$. The post-selection joint density of $(\sqrt{n}\bar{y},\omega)$ equals to
\begin{equation*}
	\frac{\phi_{(0,1)}(\sqrt{n}(\bar{y}-\mu_n))\cdot \phi_{(0,1)}(\omega/\gamma)}{(\mathbb{P}_{\sqrt{n}\bar{y}\sim\mathcal{N}(\sqrt{n}\mu_n, 1)}\times\mathbb{P}_{\omega\sim\mathcal{N}(0,\gamma^2)})\left\{\sqrt{n}\bar{y}+\omega>\tau_n\right\}}\cdot\mathbb{I}_{\{\sqrt{n}\bar{y}+\omega>\tau_n\}}.
\end{equation*}
Marginalizing over $\omega$, the post-selective density of $\sqrt{n}\bar{y}$ equals to
\begin{equation} \label{eq:simple:problem:data:density}
	\frac{\phi_{(0,1)}(\sqrt{n}(\bar{y}-\mu_n))}{1-\Phi\left(\frac{\tau_n-\sqrt{n}\mu_n}{\sqrt{1+\gamma^2}}\right)}\cdot\mathbb{P}_{\omega\sim\mathcal{N}(0,\gamma^2)}\left\{\sqrt{n}\bar{y}+\omega>\tau_n\right\}.
\end{equation}
The \textit{selective MLE} in this case is defined as the value $\hat{\mu}_n$ that maximizes the above density with respect to $\mu_n$, i.e.~$\hat{\mu}_n$ satisfies
\begin{equation*}
	\sqrt{n}(\bar{y}-\hat{\mu}_n)=\frac{\phi_{(0,1)}\left(\frac{\tau_n-\sqrt{n}\hat{\mu}_n}{\sqrt{1+\gamma^2}}\right)}{1-\Phi\left(\frac{\tau_n-\sqrt{n}\hat{\mu}_n}{\sqrt{1+\gamma^2}}\right)}\cdot\frac{1}{\sqrt{1+\gamma^2}}
\end{equation*}
by setting the derivative of the logarithm of the density in \eqref{eq:simple:problem:data:density} with respect to $\mu_n$ to zero.

We now state the result saying the selective MLE is consistent for the mean parameter in the case of Gaussian randomization and Gaussian data.

\begin{lemma}\textnormal{\textbf{(Gaussian randomization and Gaussian data: consistency of selective MLE \newline post-selection)}}
	Assuming $F_n=\mathcal{N}(\mu,1)$ for all $n$, and $G=\mathcal{N}(0,\gamma^2)$, we have
	\begin{equation*}
		\forall\epsilon>0, \;\;\;\underset{n\rightarrow\infty}{\lim} F_n^*\left\{|\hat{\mu}_n-\mu|>\epsilon\right\}=0.
	\end{equation*}
\end{lemma}

Note that the sequence of mean parameters does not change with $n$, i.e.~$\mu_n=\mu$ for all $n\geq 1$. This result is presented in \cite{selective_bayesian} for regression problems and general affine selection events. Also note that computing selective MLE requires knowing the probability of selection which has explicit form in the simple problem presented here. In more complicated selection problems, that is not the case. \cite{selective_bayesian} propose using an approximation to this selection probability. They call the resulting estimator approximate selective MLE and prove it is also consistent under the same assumptions as above.

\end{document}